\documentclass{gtmon_a}
\pdfoutput=1

\allowdisplaybreaks
\usepackage{amscd}

\proceedingstitle{Proceedings of the School and Conference in Algebraic
Topology (The Vietnam National University, Hanoi, 9-20 August 2004)}
\conferencestart{9 August 2004}
\conferenceend{20 August 2004}
\conferencename{School and Conference in Algebraic Topology}
\conferencelocation{Vietnam National University, Hanoi, Vietnam}

\editor{John Hubbuck}
\givenname{John}
\surname{Hubbuck}

\editor{Nguy\~\ecircumflex{}n H V H\uhorn{}ng}
\givenname{H\uhorn{}ng}
\surname{Nguy\~\ecircumflex{}n}

\editor{Lionel Schwartz}
\givenname{Lionel}
\surname{Schwartz}

\title[On the Stiefel--Whitney classes of the representations]{On the
Stiefel--Whitney classes\\of the representations associated with
$\mathrm{Spin}(15)$}

\author{Mamoru Mimura}
\givenname{Mamoru}
\surname{Mimura}
\address{Department of Mathematics\\
Faculty of Science\\\newline
Okayama University\\
3-1 Tsushima-naka, Okayama 700-8530\\Japan}
\email{mimura@math.okayama-u.ac.jp}
\urladdr{}

\author{Tetsu Nishimoto}
\givenname{Tetsu}
\surname{Nishimoto}
\address{Department of Welfare Business\\
Kinki Welfare University\\\newline
Fukusaki-cho, Hyogo 679-2217\\Japan}
\email{nishimoto@kinwu.ac.jp}
\urladdr{}

\volumenumber{11}
\issuenumber{}
\publicationyear{2007}
\papernumber{7}
\startpage{141}
\endpage{176}

\doi{}
\MR{}
\Zbl{}

\arxivreference{}

\keyword{Stiefel--Whitney class}
\keyword{classifying space}
\keyword{spinor group}
\keyword{exceptional Lie group}
\keyword{representation}
\subject{primary}{msc2000}{55R40}
\subject{secondary}{msc2000}{22E46}

\received{25 February 2005}
\revised{}
\accepted{17 June 2005}
\published{14 November 2007}
\publishedonline{14 November 2007}
\proposed{}
\seconded{}
\corresponding{}
\editor{}
\version{}

\makeatletter
\def\cnewtheorem#1[#2]#3{\newtheorem{#1}{#3}[section]
\expandafter\let\csname c@#1\endcsname\c@thm}

\AtBeginDocument{\let\bar\wbar\let\tilde\wtilde\let\hat\what}
\makeop{Spin}
\makeop{Sq}

\numberwithin{equation}{section}
\theoremstyle{plain}
\newtheorem{thm}{Theorem}[section]
\cnewtheorem{cor}[thm]{Corollary}
\cnewtheorem{lem}[thm]{Lemma}
\cnewtheorem{prop}[thm]{Proposition}
\cnewtheorem{claim}[thm]{Claim}
\cnewtheorem{conj}[thm]{Conjecture}
\theoremstyle{definition}
\cnewtheorem{defn}[thm]{Definition}
\cnewtheorem{notation}[thm]{Notation}
\cnewtheorem{rem}[thm]{Remark}

\makeatother  

\begin{document}

\begin{htmlabstract}
We determine the Stiefel&ndash;Whitney classes of the second exterior
representation and the spin representation of Spin(15),
which are useful to calculate the mod 2 cohomology of the classifying
space of the exceptional Lie group E<sub>8</sub>.
\end{htmlabstract}

\begin{abstract}
We determine the Stiefel--Whitney classes of the second exterior
representation and the spin representation of $\mathrm{Spin}(15)$,
which are useful to calculate the mod 2 cohomology of the classifying
space of the exceptional Lie group $E_8$.
\end{abstract}

\maketitle

\section{Introduction}
\label{sec:introduction}

The study of the cohomology of the classifying space of the Lie
groups has a long history; in particular,
among the exceptional Lie groups $G_2$, $F_4$, $E_6$, $E_7$, $E_8$,
Borel \cite{borel} first determined the algebra
structure of the mod 2 cohomology of the classifying spaces
$BG_2$ and $BF_4$ using the Serre spectral sequence.
In these cases, it is well known that the numbers of generators as
an algebra over the Steenrod algebra are 1 and 2 respectively.
Kono, Shimada and the first author \cite{KM,KMS2} studied the
mod 2 cohomology of $BE_6$ and $BE_7$, using the Rothenberg--Steenrod
spectral sequence $\{ E_r \}$ such that
\begin{displaymath}
  E_2 = \mathrm{Cotor}_A (\mathbb Z/2, \mathbb Z/2)
  \Longrightarrow H^*(BE_i;\mathbb Z/2),
\end{displaymath}
where $A = H^*(E_i;\mathbb Z/2)$ for $i = 6, 7$.
Especially, Toda \cite{toda} announced that the numbers of
generators as an algebra over the Steenrod algebra are 2 for both
 $BE_6$ and $BE_7$.
Obviously, the first generators are of degree 4 in the cases of $BG_2$, $BF_4$,
$BE_6$ and $BE_7$.
It is possible to represent the second generators as the
characteristic classes $w_{16}(\rho_4)$, $c_{16}(\rho_6)$ and
$p_{16}(\rho_7)$ of some representations
\begin{align*}
  & \rho_4 \co F_4 \longrightarrow SO(26), \\
  & \rho_6 \co E_6 \longrightarrow SU(27), \\
  & \rho_7 \co E_7 \longrightarrow Sp(28),
\end{align*}
in the case of $BF_4$, $BE_6$ and $BE_7$ respectively
(see Adams \cite[Corollaries 8.1, 8.3 and 8.2]{adams2} for the representations).

Little is known about the structure of mod 2 cohomology
of the classifying space $BE_8$ of the exceptional Lie group $E_8$.
However, it is quite natural to conjecture that the number of
generators as an algebra over the Steenrod algebra is 2.
\begin{conj}
  The mod 2 cohomology of the classifying space $BE_8$ has two
  generators as an algebra over the Steenrod algebra;
  the first one is of degree 4, and the second one is
  the 128--th Stiefel--Whitney class $w_{128}(\rho_8)$ of the adjoint
  representation
  \begin{displaymath}
    \rho_8 = \mathrm{Ad}_{E_8} \co E_8 \longrightarrow SO(248).
  \end{displaymath}
\end{conj}

(For the adjoint representation $\mathrm{Ad}_{E_8}$, see
Adams \cite[Chapters 6 and 7]{adams2}.)

Based on the computation in Mori \cite{Mori}, we conjecture more precisely
that the number of generators as an algebra (not as an algebra over
the Steenrod algebra) is 33 and the action of the Steenrod square
is expressed in the following diagram:

  
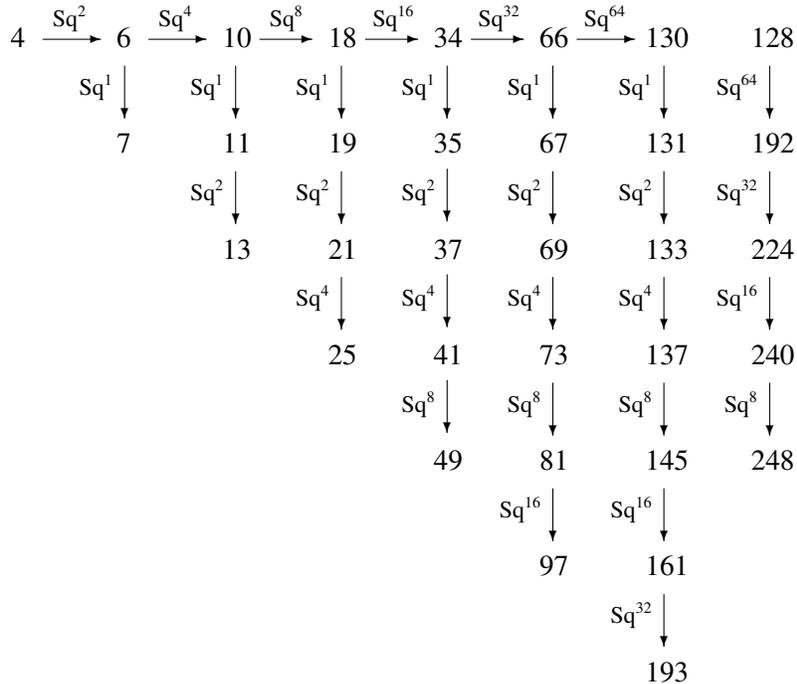
\begin{center}
    \begin{picture}(300,270)
      \put(0,250){4}
      \put(12,253){\vector(1,0){22}}
      \put(15,258){\footnotesize $\Sq^2$}

      \put(40,250){6}
      \put(52,253){\vector(1,0){22}}
      \put(55,258){\footnotesize $\Sq^4$}
      \put(43,243){\vector(0,-1){20}}
      \put(25,232){\footnotesize $\Sq^1$}

      \put(40,210){7}

      \put(80,250){10}
      \put(94,253){\vector(1,0){20}}
      \put(97,258){\footnotesize $\Sq^8$}
      \put(85,243){\vector(0,-1){20}}
      \put(67,232){\footnotesize $\Sq^1$}

      \put(80,210){11}
      \put(85,203){\vector(0,-1){20}}
      \put(67,192){\footnotesize $\Sq^2$}

      \put(80,170){13}

      \put(120,250){18}
      \put(134,253){\vector(1,0){20}}
      \put(136,258){\footnotesize $\Sq^{16}$}
      \put(125,243){\vector(0,-1){20}}
      \put(107,232){\footnotesize $\Sq^1$}

      \put(120,210){19}
      \put(125,203){\vector(0,-1){20}}
      \put(107,192){\footnotesize $\Sq^2$}

      \put(120,170){21}
      \put(125,163){\vector(0,-1){20}}
      \put(107,152){\footnotesize $\Sq^4$}

      \put(120,130){25}

      \put(160,250){34}
      \put(174,253){\vector(1,0){20}}
      \put(176,258){\footnotesize $\Sq^{32}$}
      \put(165,243){\vector(0,-1){20}}
      \put(147,232){\footnotesize $\Sq^1$}

      \put(160,210){35}
      \put(165,204){\vector(0,-1){20}}
      \put(147,192){\footnotesize $\Sq^2$}

      \put(160,170){37}
      \put(165,164){\vector(0,-1){20}}
      \put(147,152){\footnotesize $\Sq^4$}

      \put(160,130){41}
      \put(165,124){\vector(0,-1){20}}
      \put(147,112){\footnotesize $\Sq^8$}

      \put(160,90){49}

      \put(200,250){66}
      \put(214,253){\vector(1,0){22}}
      \put(216,258){\footnotesize $\Sq^{64}$}
      \put(205,243){\vector(0,-1){20}}
      \put(187,232){\footnotesize $\Sq^1$}

      \put(200,210){67}
      \put(205,203){\vector(0,-1){20}}
      \put(187,192){\footnotesize $\Sq^2$}

      \put(200,170){69}
      \put(205,163){\vector(0,-1){20}}
      \put(187,152){\footnotesize $\Sq^4$}

      \put(200,130){73}
      \put(205,123){\vector(0,-1){20}}
      \put(187,112){\footnotesize $\Sq^8$}

      \put(200,90){81}
      \put(205,83){\vector(0,-1){20}}
      \put(184,72){\footnotesize $\Sq^{16}$}

      \put(200,50){97}

      \put(240,250){130}
      \put(247,243){\vector(0,-1){20}}
      \put(229,232){\footnotesize $\Sq^1$}

      \put(240,210){131}
      \put(247,203){\vector(0,-1){20}}
      \put(229,192){\footnotesize $\Sq^2$}

      \put(240,170){133}
      \put(247,163){\vector(0,-1){20}}
      \put(229,152){\footnotesize $\Sq^4$}

      \put(240,130){137}
      \put(247,123){\vector(0,-1){20}}
      \put(229,112){\footnotesize $\Sq^8$}

      \put(240,90){145}
      \put(247,83){\vector(0,-1){20}}
      \put(226,72){\footnotesize $\Sq^{16}$}

      \put(240,50){161}
      \put(247,43){\vector(0,-1){20}}
      \put(226,32){\footnotesize $\Sq^{32}$}

      \put(240,10){193}

      \put(280,250){128}
      \put(287,243){\vector(0,-1){20}}
      \put(266,232){\footnotesize $\Sq^{64}$}

      \put(280,210){192}
      \put(287,203){\vector(0,-1){20}}
      \put(266,192){\footnotesize $\Sq^{32}$}

      \put(280,170){224}
      \put(287,163){\vector(0,-1){20}}
      \put(266,152){\footnotesize $\Sq^{16}$}

      \put(280,130){240}
      \put(287,123){\vector(0,-1){20}}
      \put(269,112){\footnotesize $\Sq^8$}

      \put(280,90){248}
    \end{picture}
  \end{center}
where the numbers in the diagram indicate the degrees of the
generators.

Our main tool to prove the conjecture on $H^*(BE_8;\mathbb Z/2)$
is the Rothenberg--Steenrod spectral sequence.

Our tactics to prove the conjecture may be stated as follows; firstly we
calculate $E_2 = \mathrm{Cotor}_A (\mathbb Z/2, \mathbb Z/2)$ for $A = H^*(E_8;\mathbb 
Z/2)$.
We claim that generators are exactly given by the elements of degree
indicated in the above diagram;
we will determine it in the forthcoming paper (see Mori \cite{Mori}), and
secondly we show that these generators are detected by the images of
the Steenrod squares of the generator of degree 4 and by the images of
the Steenrod squares of the generator of degree 128.
Observe that the remaining generators except those in the right column
in the above diagram come from the mod 2 cohomology of the
Eilenberg--Mac Lane space $K(\mathbb Z, 4)$.

In order to prove the conjecture on the mod 2 cohomology of $BE_8$,
we need to find a Lie group $G$ and a homomorphism to $E_8$ such that
the structure of $H^*(BG;\mathbb Z/2)$ as an algebra over the Steenrod
algebra is
known and that the induced homomorphism $H^n(BE_8;\mathbb Z/2) \to
H^n(BG;\mathbb Z/2)$ is monic for $n \leq N$, where $N$ is sufficiently
large.
Using the homomorphism $G \to E_8$ with the above properties, we
consider the Stiefel--Whitney class of the induced representation
\begin{displaymath}
  G \longrightarrow E_8 \xrightarrow{\mathrm{Ad}_{E_8}} SO(248).
\end{displaymath}
The natural inclusion map of the semi-spinor group $Ss(16) \subset E_8$
might seem to be the best homomorphism among others, since the homomorphism
\begin{displaymath}
  H^*(BE_8;\mathbb Z/2) \longrightarrow H^*(BSs(16);\mathbb Z/2)
\end{displaymath}
is a monomorphism.
However, calculating $H^*(BSs(16);\mathbb Z/2)$ seems to be as difficult as
calculating $H^*(BE_8;\mathbb Z/2)$, since the Hopf algebra structure
of $H^*(Ss(16);\mathbb Z/2)$ is similar to that of $H^*(E_8;\mathbb Z/2)$.
The next candidate is the spinor group $\Spin(16)$ which is the
universal covering of $Ss(16)$ with the obvious homomorphism
$$\Spin(16)\rightarrow Ss(16) \rightarrow E_8.$$
According to Adams \cite{adams2}, the Lie algebra $L(E_8)$ of type $E_8$
can be constructed as the direct sum $L(\Spin(16)) \oplus \smash{\Delta_{16}^+}$
with some Lie algebra structure, where $L(\Spin(16))$ is the Lie
algebra of type $\Spin(16)$ and $\smash{\Delta_{16}^+}$ is the spin
representation of $\Spin(16)$:
\begin{displaymath}
  \Delta_{16}^+ \co \Spin(16) \longrightarrow SO(128).
\end{displaymath}
Note that the Lie algebra $L(\Spin(n))$ of type $\Spin(n)$ is isomorphic
to the second exterior representation
\begin{displaymath}
  \lambda_n^2 \co \Spin(n) \longrightarrow SO(n) \longrightarrow
  SO\big(\tbinom{n}{2}\big)
\end{displaymath}
as a representation of $\Spin(n)$.
Thus the induced representation
$\Spin(16) \rightarrow Ss(16) \rightarrow E_8 \rightarrow SO(248)$
is the direct sum $\lambda_{16}^2 \oplus \Delta_{16}^+$.
The more appropriate homomorphism than $\Spin(16) \to E_8$ is the
composition map of the natural maps
\begin{displaymath}
  \Spin(15) \longrightarrow \Spin(16) \longrightarrow Ss(16)
  \longrightarrow E_8,
\end{displaymath}
since the image of $H^*(BSs(16);\mathbb Z/2) \to H^*(B\Spin(16);
\mathbb Z/2)$ is isomorphic to that of $H^*(BSs(16);\mathbb Z/2)
\to H^*(B\Spin(15);\mathbb Z/2)$ and since $\Spin(15)$ is smaller
than $\Spin(16)$.
Observe that the induced representation
\begin{displaymath}
  \Spin(15) \longrightarrow \Spin(16) \longrightarrow Ss(16) \longrightarrow
  E_8 \xrightarrow{\mathrm{Ad}_{E_8}} SO(248)
\end{displaymath}
is a direct sum of the first exterior representation (or the projection
map)
\begin{displaymath}
  \lambda_{15}^1 \co \Spin(15) \longrightarrow SO(15),
\end{displaymath}
the second exterior representation
\begin{displaymath}
  \lambda_{15}^2 \co \Spin(15) \longrightarrow SO(105),
\end{displaymath}
and the spin representation
\begin{displaymath}
  \Delta_{15} \co \Spin(15) \longrightarrow SO(128),
\end{displaymath}
since we have $f_{15}^*\lambda_{16}^2 = \lambda_{15}^1 \oplus
\lambda_{15}^2$ and $f_{15}^*\Delta_{16}^+ = \Delta_{15}$,
where $f_{15} \co \Spin(15) \to \Spin(16)$ is the natural inclusion map.
Quillen's theorem (cf \fullref{thm:quillen}) states that the 128--th
Stiefel--Whitney class $\smash{w_{128}(\Delta_{15})}$ of the spin
representation $\smash{\Delta_{15}}$ is a member of a system
of generators of $\smash{H^*(B\Spin(15);\mathbb Z/2)}$ as an algebra over
the Steenrod algebra.
Our calculation asserts (\fullref{thm:adjoint}) that
$$w_{128}(\lambda_{15}^1 \oplus
\lambda_{15}^2 \oplus \Delta_{15}) \equiv w_{128}(\Delta_{15}) \text{ mod
decomposables,}$$
which implies (\fullref{cor:gen}) that the 128--th Stiefel--Whitney class
$w_{128}(\mathrm{Ad}_{E_8})$ can be chosen as a member of a system
of generators of $H^*(BE_8;\mathbb Z/2)$ as an algebra over the Steenrod
algebra.

The paper is organized as follows.
In \fullref{sec:preliminary}, we prepare some results on representations which will
be needed for later use.
In \fullref{sec:2nd_exterior_rep}, we calculate the
Stiefel--Whitney classes of the second exterior representation
$\lambda_{15}^2$ of $\Spin(15)$.
In \fullref{sec:spin_rep}, we calculate the Stiefel--Whitney
classes of the spin representation $\Delta_{15}$ of $\Spin(15)$.
One can derive the total characteristic classes $w(\rho_4)$, $c(\rho_6)$,
$p(\rho_7)$ from the results \eqref{eq:delta9}, \eqref{eq:delta10},
\eqref{eq:delta12-} respectively.
In \fullref{sec:adjoint_rep}, we calculate the Stiefel--Whitney
classes of the representation $\lambda_{15}^1 \oplus \lambda_{15}^2
\oplus \Delta_{15}$ of $\Spin(15)$ which is the induced representation
from the adjoint representation of $E_8$.

Our main results are \fullref{thm:adjoint} and \fullref{cor:gen};
the latter result assures the existence of the first generator in the
right column in the previous diagram. The former result will be
used to show in the forthcoming paper that the action of the
Steenrod squares on this generator of degree 128, and hence on those
of degrees 192, 224, 240, 248, by using the Wu formula. That is,
the generators in the right column
can be represented by the 128--th, 192--nd, 224--th, 240--th and
248--th Stiefel--Whitney classes of the adjoint representation
respectively.

It is our pleasure to acknowledge that the present paper is motivated
by the calculation in Mori \cite{Mori}.
Most of the calculations were performed by programs using GAP which is a
system for computation in discrete abstract algebra.
We thank Shingo Okuyama and Yuriko Sambe who advised us about
programming of the calculations.

\section{Preliminary}
\label{sec:preliminary}

In this section, we recall the mod 2 cohomology of the classifying
spaces of $O(n)$, $SO(n)$ and $\Spin(n)$ as well as the
Stiefel--Whitney classes of representations.

Let $H_n$ be the subgroup of $O(n)$ consisting of the diagonal
matrices, which is isomorphic to $(\mathbb Z/2)^n$.
Let $i_n \co H_n \to O(n)$ be the natural inclusion map.
Let $W = N(H_n)/H_n$ be the Weyl group of $O(n)$, where $N(H_n)$ is
the normalizer of $H_n$.
As is well known, $W$ is isomorphic to the $n$--th symmetric
group $\Sigma_n$.
The mod 2 cohomology of $BO(n)$ is a polynomial algebra whose
generators are defined as the invariants under the action of
the Weyl group:
\begin{align*}
  H^*(BO(n);\mathbb Z/2) & = H^*(B(\mathbb Z/2)^n;\mathbb Z/2)^W \\
  & = \smash{\mathbb Z/2[t_1, t_2, \ldots, t_n]^W} \\
  & = \mathbb Z/2[w_1, w_2, \ldots, w_n],
\end{align*}
where $\{t_j : 1 \leq j \leq n\}$ is a basis of $H^1(B(\mathbb
Z/2)^n;\mathbb Z/2)$ and $w_i$, the $i$--th elementary symmetric
polynomial of $t_j$, is called the $i$--th Stiefel--Whitney class.
Similarly, the mod 2 cohomology of $BSO(n)$ is a polynomial algebra
generated by $w_i$ for $2 \leq i \leq n$:
\begin{displaymath}
  H^*(BSO(n);\mathbb Z/2) = \mathbb Z/2[w_2, w_3, \ldots, w_n].
\end{displaymath}
Let $G$ be a compact Lie group and $\rho$ an $n$--dimensional real
representation $G \to O(n)$.
Since a homomorphism induces a map between their classifying spaces
\begin{displaymath}
  B\rho \co BG \longrightarrow BO(n)
\end{displaymath}
uniquely up to homotopy, we obtain a homomorphism between their mod 2
cohomologies:
\begin{displaymath}
  B\rho^* \co H^*(BO(n);\mathbb Z/2) \longrightarrow H^*(BG;\mathbb Z/2).
\end{displaymath}
We denote $B\rho^*(w_i)$ simply by $w_i(\rho)$.
For a representation $\rho \co G \to SO(n)$, we also denote
by $w_i(\rho)$ the induced element.
One of the important properties of the Stiefel--Whitney class is
the Whitney product formula:
\begin{equation}
  \label{eq:whitney_sum}
  w_k(\iota_{m,n}) = \sum_{i+j=k} w'_i \times w''_j,
\end{equation}
where $\iota_{m,n} \co O(m) \times O(n) \to O(m+n)$ is the obvious map,
and $w'_i$ and
$w''_i$ are the $i$--th Stiefel--Whitney classes of $H^*(BO(m);\mathbb
Z/2)$ and $H^*(BO(n);\mathbb Z/2)$ respectively.

The action of the Steenrod square on the
Stiefel--Whitney classes is given by the Wu formula:
\begin{displaymath}
  \Sq^j w_i = \sum_{k=0}^j \binom{i-k-1}{j-k} w_{i+j-k}w_k \quad (0
  \leq j \leq i).
\end{displaymath}
Then, using it, one can easily see that generators of
$H^*(BO(n);\mathbb Z/2)$ as an algebra over the Steenrod algebra
are given by $w_{2^k}$ for $1 \leq 2^k \leq n$.

Now we recall a result due to Borel--Hirzebruch \cite{BH1}.
Let $H$ be an elementary abelian 2--subgroup of $G$, and $i \co H \to G$
the inclusion map.
Suppose that there is a representation $\rho \co G \to O(n)$ satisfying
the following commutative diagram:
\begin{displaymath}
  \begin{CD}
    H @>{\bar{\rho}}>> H_n \\
    @V{i}VV @VV{i_n}V \\
    G @>>{\rho}> O(n),
  \end{CD}
\end{displaymath}
where $\bar{\rho} = \rho|_H$.
\begin{prop}
  \label{prop:BH}
  {\rm (Borel--Hirzebruch \cite{BH1})}
  There holds
  \begin{displaymath}
    Bi^*(w(\rho)) = \prod_{i=1}^n (1+B\bar\rho^*(t_i)).
  \end{displaymath}
  Moreover, if $G = O(m)$ and $H = H_m$, then $Bi^*(w(\rho))$ is a
  symmetric polynomial of $t'_i$ for $1 \leq i \leq m$, where $t'_i \in
  H^1(BH_m;\mathbb Z/2)$.
\end{prop}
Thus one can calculate $w(\rho)$,
if $B\bar\rho^*(t_i)$ is calculable.

To state Quillen's result \cite{quillen} concerning
$H^*(B\Spin(n);\mathbb Z/2)$,
we define $h_n$ as follows:
\begin{displaymath}
  h_n =
  \begin{cases}
    (n-2)/2, & n \equiv 0 \mod 8, \\
    (n-1)/2, & n \equiv 1, 7 \mod 8, \\
    n/2, & n \equiv 2, 4, 6 \mod 8, \\
    (n+1)/2, & n \equiv 3, 5 \mod 8.
  \end{cases}
\end{displaymath}
We often denote $h_n$ by $h$ if there is no confusion.
Let $J \subset H^*(BSO(n);\mathbb Z/2)$ be the ideal generated by
the elements
\begin{displaymath}
  w_2, \hspace{7pt} \Sq^1 w_2, \hspace{7pt} \Sq^2\Sq^1 w_2, \hspace{7pt}
  \ldots,
  \hspace{7pt} \Sq^{2^{h-2}}\Sq^{2^{h-3}} \cdots \Sq^2\Sq^1 w_2.
\end{displaymath}
As is well known (see for example Adams \cite{adams2}),
$\Spin(n)$ has the spin representations:
\begin{align*}
  & \Delta_{8m}^{\pm} \co \Spin(8m) \longrightarrow SO(2^{4m-1}), \\
  & \Delta_{8m+1} \co \Spin(8m+1) \longrightarrow SO(2^{4m}), \\
  & \Delta_{8m+2}^{\pm} \co \Spin(8m+2) \longrightarrow SU(2^{4m}), \\
  & \Delta_{8m+3} \co \Spin(8m+3) \longrightarrow Sp(2^{4m}), \\
  & \Delta_{8m+4}^{\pm} \co \Spin(8m+4) \longrightarrow Sp(2^{4m}), \\
  & \Delta_{8m+5} \co \Spin(8m+5) \longrightarrow Sp(2^{4m+1}), \\
  & \Delta_{8m+6}^{\pm} \co \Spin(8m+6) \longrightarrow SU(2^{4m+2}), \\
  & \Delta_{8m+7} \co \Spin(8m+7) \longrightarrow SO(2^{4m+3}),
\end{align*}
where $\Delta^{\pm}$ means that there are two representations
$\Delta^+$ and $\Delta^-$.
\begin{notation}
  \begin{enumerate}
  \item For a real representation $\rho \co G \to SO(n)$, we denote by
    $\rho_{\mathbb C}$ a complex representation
    \begin{displaymath}
      G \overset{\rho}{\longrightarrow} SO(n) \longrightarrow SU(n).
    \end{displaymath}
  \item For a quaternionic representation $\rho \co G \to Sp(n)$, we
    denote by $\rho_{\mathbb C}$ a complex representation
    \begin{displaymath}
      G \overset{\rho}{\longrightarrow} Sp(n) \longrightarrow SU(2n).
    \end{displaymath}
  \item For a complex representation $\rho \co G \to SU(n)$, we denote by
    $\rho_{\mathbb R}$ a real representation
    \begin{displaymath}
      G \overset{\rho}{\longrightarrow} SU(n) \longrightarrow SO(2n).
    \end{displaymath}
  \item For a quaternionic representation $\rho \co G \to Sp(n)$, we
    denote by $\rho_{\mathbb R}$ a real representation
    \begin{displaymath}
      G \overset{\rho}{\longrightarrow} Sp(n) \longrightarrow SU(2n)
      \longrightarrow SO(4n).
    \end{displaymath}
  \end{enumerate}
\end{notation}
The following remark is well known (see for example Adams \cite{adams1}).
\begin{rem}
  \label{rem:iso}
  For real representations $\rho$ and $\sigma$,
  $\rho_{\mathbb C} \cong \sigma_{\mathbb C}$ if and only if
  $\rho \cong \sigma$.
  For quaternionic representations $\rho$ and $\sigma$,
  $\rho_{\mathbb C} \cong \sigma_{\mathbb C}$ if and only if
  $\rho \cong \sigma$.
\end{rem}
We sometimes denote by $\rho_{\mathbb R}$ a real representation
$\rho$, and denote by $\rho_{\mathbb C}$ a complex representation
$\rho$, by abuse of the notations.
It is well known that $w(\rho_{\mathbb R}) = c(\rho)$ if $\rho$ is a
complex representation.
It is also known that $w(\rho_{\mathbb R}) = c(\rho_{\mathbb C}) =
p(\rho)$ if $\rho$ is a quaternionic representation.
Note that, if $n = 4m+2$, the representations of $(\Delta_n^+)_{\mathbb
  R}$ and $(\Delta_n^-)_{\mathbb R}$ are isomorphic, since
$\Delta_n^+$ and $\Delta_n^-$ are conjugate to each other.
We denote $w_{2^h}((\Delta_n)_{\mathbb R})$ and
$w_{2^h}((\Delta_n^+)_{\mathbb R})$ simply by
$u_{2^h}$ for $n$ odd and $n$ even respectively.

Now the results due to Quillen \cite{quillen} can be summarized as
follows:
\begin{thm}
  \label{thm:quillen}
  {\rm (1)}\qua
  The algebra structure of the {\rm mod 2} cohomology of $B\Spin(n)$ is
  given by
  \begin{displaymath}
    H^*(B\Spin(n);\mathbb Z/2) \cong
    H^*(BSO(n);\mathbb Z/2)/J \otimes
    \mathbb Z/2[u_{2^h}].
  \end{displaymath}
  {\rm (2)}\qua 
  The nonzero Stiefel--Whitney classes of the spin representation are
  those of degrees $2^h$ and $2^h-2^i$ for $r \leq i \leq h$,  where
  \begin{displaymath}
    r =
    \begin{cases}
      0 & n \equiv 0, 1, 7 \mod 8, \\
      1 & n \equiv 2, 6\phantom{, 6}  \mod 8, \\
      2 & n \equiv 3, 4, 5 \mod 8.
    \end{cases}
  \end{displaymath}
\end{thm}

In order to choose a Gr\"obner basis of the ideal of $J$ (see for
example Cox, Little and O'Shea \cite{CLO}), we need to introduce a total order to the set of monomials in
$H^*(BSO(n);\mathbb Z/2)$.
Let $a = w_2^{j_2}w_3^{j_3}\cdots w_n^{j_n}$ and $b =
w_2^{k_2}w_3^{k_3}\cdots w_n^{k_n}$ be monomials in
$H^*(BSO(n);\mathbb Z/2)$.
Then we define a total order as follows:
\begin{displaymath}
  a < b \Longleftrightarrow
  \begin{cases}
    \deg a < \deg b, \textrm{ or} \\
    \deg a = \deg b,
    \  j_2 = k_2,
    \  \ldots,
    \  j_{i-1} = k_{i-1},
    \  j_i > k_i.
  \end{cases}
\end{displaymath}
For the projection $\lambda^1 \co \Spin(n) \to SO(n)$,
denote $w_i(\lambda^1)$ simply by $y_i$, since we need to distinguish
the elements of $H^*(BSO(n);\mathbb Z/2)$ from those of
$H^*(B\Spin(n);\mathbb Z/2)$.
We exclude the elements $y_2$, $y_3$, $y_5$ and $y_9$ in the generators of
$H^*(B\Spin(n);\mathbb Z/2)$, since $\Sq^1 y_2 = y_3$, $\Sq^2 y_3 = y_5 +
y_2y_3$ and $\Sq^4 y_5 = y_9 + y_2y_7 + y_3y_6 + y_4y_5$,
that is, they are 0 in $H^*(B\Spin(n);\mathbb Z/2)$.
Now we can describe the algebra structure of $H^*(B\Spin(n);\mathbb Z/2)$ for
$3 \leq n \leq 15$ explicitly using the reduced Gr\"obner basis $\{ R_i \}$
as follows:
\begin{align*}
  & H^*(B\Spin(3);\mathbb Z/2) = \mathbb Z/2[u_4]; \\
  & H^*(B\Spin(4);\mathbb Z/2) = \mathbb Z/2[y_4, u_4]; \\
  & H^*(B\Spin(5);\mathbb Z/2) = \mathbb Z/2[y_4, u_8]; \\
  & H^*(B\Spin(6);\mathbb Z/2) = \mathbb Z/2[y_4, y_6, u_8]; \\
  & H^*(B\Spin(7);\mathbb Z/2) = \mathbb Z/2[y_4, y_6, y_7, u_8]; \\
  & H^*(B\Spin(8);\mathbb Z/2) = \mathbb Z/2[y_4, y_6, y_7, y_8, u_8]; \\
  & H^*(B\Spin(9);\mathbb Z/2) = \mathbb Z/2[y_4, y_6, y_7, y_8, u_{16}]; \\
  & H^*(B\Spin(10);\mathbb Z/2) = \mathbb Z/2[y_4, y_6, y_7, y_8,
  y_{10}, u_{32}]/(R_1), \\
  & \ \
  R_1 = y_7y_{10}; \\
  & H^*(B\Spin(11);\mathbb Z/2) = \mathbb Z/2[y_4, y_6, y_7, y_8,
  y_{10}, y_{11}, u_{64}]/(R_1, R_2), \\
  & \ \ 
  R_1 = y_7y_{10}{+}y_6y_{11}, \\
  & \ \ 
  R_2 = y_{11}^3{+}y_7^2y_8y_{11}{+}y_4y_7y_{11}^2; \\
  & H^*(B\Spin(12);\mathbb Z/2) = \mathbb Z/2[y_4, y_6, y_7, y_8,
  y_{10}, y_{11}, y_{12}, u_{64}]/(R_1, R_2), \\
  & \ \ 
  R_1 = y_7y_{10}{+}y_6y_{11}, \\
  & \ \ 
  R_2 = y_{11}^3{+}y_7^2y_8y_{11}{+}y_7^3y_{12}{+}y_4y_7y_{11}^2; \\
  & H^*(B\Spin(13);\mathbb Z/2) = \mathbb Z/2[y_4, y_6, y_7, y_8,
  y_{10}, y_{11}, y_{12}, y_{13}, u_{128}]/(R_1, R_2, R_3), \\
  & \ \ 
  R_1 =
  y_7y_{10}
 {+}y_6y_{11}
 {+}y_4y_{13}
  , \\
  & \ \ 
  R_2 =
  y_{11}^3
 {+}y_{10}^2y_{13}
 {+}y_7y_{13}^2
 {+}y_7^2y_8y_{11}
 {+}y_7^3y_{12}
 {+}y_6y_7^2y_{13}
 {+}y_6^2y_8y_{13}
 {+}y_4y_7y_{11}^2
  \\ & \qquad
 {+}y_4y_6y_{10}y_{13}
 {+}y_4^2y_{12}y_{13}
  , \\
  & \ \ 
  R_3 =
  y_{13}^5
 {+}y_{10}^3y_{11}^2y_{13}
 {+}y_{10}^4y_{12}y_{13}
 {+}y_8y_{10}^2y_{11}y_{13}^2
 {+}y_7y_8y_{11}y_{13}^3
 {+}y_7^2y_8^2y_{11}^2y_{13}
  \\ & \qquad
 {+}y_7^3y_8y_{11}y_{12}y_{13}
 {+}y_7^4y_{12}^2y_{13}
 {+}y_7^4y_8^3y_{13}
 {+}y_6y_{10}^2y_{13}^3
 {+}y_6y_7^2y_8y_{11}y_{13}^2
  \\ & \qquad
 {+}y_6^2y_8^2y_{11}y_{13}^2
 {+}y_6^2y_7^2y_{13}^3
 {+}y_6^3y_{11}^2y_{12}y_{13}
 {+}y_6^3y_{10}y_{11}y_{13}^2
 {+}y_6^3y_8y_{13}^3
  \\ & \qquad
 {+}y_6^3y_7y_8^2y_{11}y_{13}
 {+}y_6^3y_7^2y_8y_{12}y_{13}
 {+}y_6^4y_8^2y_{12}y_{13}
 {+}y_6^4y_7y_8y_{13}^2
 {+}y_4y_{11}^2y_{13}^3
  \\ & \qquad
 {+}y_4y_7^2y_{11}^2y_{12}y_{13}
 {+}y_4y_7^3y_8^2y_{11}y_{13}
 {+}y_4y_6y_7y_{11}^2y_{13}^2
 {+}y_4y_6^3y_7y_{11}y_{12}y_{13}
  \\ & \qquad
 {+}y_4y_6^4y_{11}y_{13}^2
 {+}y_4^2y_{10}y_{11}^2y_{12}y_{13}
 {+}y_4^2y_{10}^2y_{11}y_{13}^2
 {+}y_4^2y_8y_{11}y_{12}y_{13}^2
  \\ & \qquad
 {+}y_4^2y_7y_{11}y_{13}^3
 {+}y_4^2y_7^3y_{11}y_{12}y_{13}
 {+}y_4^2y_6y_7y_8y_{11}y_{12}y_{13}
 {+}y_4^2y_6y_7^2y_{12}^2y_{13}
  \\ & \qquad
 {+}y_4^2y_6y_7^2y_{11}y_{13}^2
 {+}y_4^2y_6^2y_{10}y_{11}^2y_{13}
 {+}y_4^2y_6^2y_{10}^2y_{12}y_{13}
 {+}y_4^2y_6^3y_{13}^3
  \\ & \qquad
 {+}y_4^3y_6y_{11}^2y_{12}y_{13}
 {+}y_4^3y_6y_{10}y_{11}y_{13}^2
 {+}y_4^4y_{11}y_{12}y_{13}^2
 {+}y_4^4y_{12}^3y_{13}
  ; \\
  & H^*(B\Spin(14);\mathbb Z/2) = \\
  &\hphantom{H^*(B\Spin(14);}\mathbb Z/2[y_4, y_6, y_7, y_8,
  y_{10}, y_{11}, y_{12}, y_{13}, y_{14}, u_{128}]/(R_1, R_2, R_3), \\
  & \ \ 
  R_1 =
  y_7y_{10}
 {+}y_6y_{11}
 {+}y_4y_{13}
  , \\
  & \ \ 
  R_2 =
  y_{11}^3
 {+}y_{10}^2y_{13}
 {+}y_7y_{13}^2
 {+}y_7^2y_8y_{11}
 {+}y_7^3y_{12}
 {+}y_6y_7^2y_{13}
 {+}y_6^2y_8y_{13}
 {+}y_6^2y_7y_{14}
  \\ & \qquad
 {+}y_4y_7y_{11}^2
 {+}y_4y_6y_{10}y_{13}
 {+}y_4^2y_{12}y_{13}
 {+}y_4^2y_{11}y_{14}
  , \\
  & \ \ 
  R_3 =
  y_{13}^5
 {+}y_{10}^3y_{11}^2y_{13}
 {+}y_{10}^4y_{12}y_{13}
 {+}y_{10}^4y_{11}y_{14}
 {+}y_8y_{10}^2y_{11}y_{13}^2
 {+}y_7y_8y_{11}y_{13}^3
  \\ & \qquad
 {+}y_7^2y_{11}y_{13}^2y_{14}
 {+}y_7^2y_8^2y_{11}^2y_{13}
 {+}y_7^3y_8y_{11}y_{12}y_{13}
 {+}y_7^3y_8y_{11}^2y_{14}
 {+}y_7^4y_{12}^2y_{13}
  \\ & \qquad
 {+}y_7^4y_{11}y_{12}y_{14}
 {+}y_7^4y_8^3y_{13}
 {+}y_7^5y_8^2y_{14}
 {+}y_6y_{10}y_{11}^2y_{13}y_{14}
 {+}y_6y_{10}^2y_{13}^3
  \\ & \qquad
 {+}y_6y_7^2y_8y_{11}y_{13}^2
 {+}y_6y_7^3y_{11}y_{13}y_{14}
 {+}y_6^2y_8^2y_{11}y_{13}^2
 {+}y_6^2y_7^2y_{13}^3
 {+}y_6^2y_7^2y_{11}y_{14}^2
  \\ & \qquad
 {+}y_6^3y_{11}^2y_{12}y_{13}
 {+}y_6^3y_{10}y_{11}y_{13}^2
 {+}y_6^3y_{10}^2y_{13}y_{14}
 {+}y_6^3y_8y_{13}^3
 {+}y_6^3y_7y_8^2y_{11}y_{13}
  \\ & \qquad
 {+}y_6^3y_7^2y_8y_{12}y_{13}
 {+}y_6^4y_{13}y_{14}^2
 {+}y_6^4y_8^2y_{12}y_{13}
 {+}y_6^4y_8^2y_{11}y_{14}
 {+}y_6^4y_7y_8y_{13}^2
  \\ & \qquad
 {+}y_6^5y_8y_{13}y_{14}
 {+}y_4y_{11}^2y_{13}^3
 {+}y_4y_{10}y_{11}y_{13}^2y_{14}
 {+}y_4y_7^2y_{11}^2y_{12}y_{13}
 {+}y_4y_7^3y_{13}^2y_{14}
  \\ & \qquad
 {+}y_4y_7^3y_8^2y_{11}y_{13}
 {+}y_4y_7^4y_8y_{11}y_{14}
 {+}y_4y_7^5y_{12}y_{14}
 {+}y_4y_6y_7y_{11}^2y_{13}^2
  \\ & \qquad
 {+}y_4y_6y_7^4y_{13}y_{14}
 {+}y_4y_6^2y_{11}^2y_{13}y_{14}
 {+}y_4y_6^2y_7^2y_8y_{13}y_{14}
 {+}y_4y_6^2y_7^3y_{14}^2
  \\ & \qquad
 {+}y_4y_6^3y_7y_{11}y_{12}y_{13}
 {+}y_4y_6^4y_{11}y_{13}^2
 {+}y_4y_6^4y_{10}y_{13}y_{14}
 {+}y_4^2y_{10}y_{11}^2y_{12}y_{13}
  \\ & \qquad
 {+}y_4^2y_{10}^2y_{11}y_{13}^2
 {+}y_4^2y_{10}^3y_{13}y_{14}
 {+}y_4^2y_8y_{11}y_{12}y_{13}^2
 {+}y_4^2y_8y_{11}^2y_{13}y_{14}
  \\ & \qquad
 {+}y_4^2y_7y_{11}y_{13}^3
 {+}y_4^2y_7y_{11}y_{12}y_{13}y_{14}
 {+}y_4^2y_7y_{11}^2y_{14}^2
 {+}y_4^2y_7^3y_{11}y_{12}y_{13}
  \\ & \qquad
 {+}y_4^2y_6y_{11}y_{13}^2y_{14}
 {+}y_4^2y_6y_7y_8y_{11}y_{12}y_{13}
 {+}y_4^2y_6y_7^2y_{12}^2y_{13}
 {+}y_4^2y_6y_7^2y_{11}y_{13}^2
  \\ & \qquad
 {+}y_4^2y_6^2y_{10}y_{11}^2y_{13}
 {+}y_4^2y_6^2y_{10}^2y_{12}y_{13}
 {+}y_4^2y_6^2y_{10}^2y_{11}y_{14}
 {+}y_4^2y_6^2y_8y_{10}y_{13}y_{14}
  \\ & \qquad
 {+}y_4^2y_6^3y_{12}y_{13}y_{14}
 {+}y_4^2y_6^3y_{13}^3
 {+}y_4^3y_6y_{11}^2y_{12}y_{13}
 {+}y_4^3y_6y_{10}y_{11}y_{13}^2
  \\ & \qquad
 {+}y_4^3y_6y_{10}^2y_{13}y_{14}
 {+}y_4^3y_6^2y_{13}y_{14}^2
 {+}y_4^3y_7y_8y_{11}y_{13}y_{14}
 {+}y_4^3y_7^2y_{11}y_{14}^2
  \\ & \qquad
 {+}y_4^4y_{10}y_{12}y_{13}y_{14}
 {+}y_4^4y_{11}y_{12}y_{13}^2
 {+}y_4^4y_{11}y_{12}^2y_{14}
 {+}y_4^4y_{12}^3y_{13}
 {+}y_4^4y_8y_{13}y_{14}^2
  \\ & \qquad
 {+}y_4^4y_7y_{14}^3
  ; \\
  & H^*(B\Spin(15);\mathbb Z/2) = \\
  &\hphantom{H^*(B\Spin(15);}\mathbb Z/2[y_4, y_6, y_7, y_8,
  y_{10}, y_{11}, y_{12}, y_{13}, y_{14}, y_{15}, u_{128}]/(R_1, R_2, R_3), \\
  & \ \ 
  R_1 =
  y_7y_{10}
 {+}y_6y_{11}
 {+}y_4y_{13}
  , \\
  & \ \ 
  R_2 =
  y_{11}^3
 {+}y_{10}^2y_{13}
 {+}y_7y_{13}^2
 {+}y_7^2y_8y_{11}
 {+}y_7^3y_{12}
 {+}y_6y_7^2y_{13}
 {+}y_6^2y_8y_{13}
 {+}y_6^2y_7y_{14}
  \\ & \qquad
 {+}y_6^3y_{15}
 {+}y_4y_7y_{11}^2
 {+}y_4y_7^2y_{15}
 {+}y_4y_6y_{10}y_{13}
 {+}y_4^2y_{12}y_{13}
 {+}y_4^2y_{11}y_{14}
  \\ & \qquad
 {+}y_4^2y_{10}y_{15}
  , \\
  & \ \ 
  R_3 =
  y_{13}^5
 {+}y_{10}^3y_{11}^2y_{13}
 {+}y_{10}^4y_{12}y_{13}
 {+}y_{10}^4y_{11}y_{14}
 {+}y_{10}^5y_{15}
 {+}y_8y_{10}^2y_{11}y_{13}^2
  \\ & \qquad
 {+}y_7y_8y_{11}y_{13}^3
 {+}y_7^2y_{11}y_{13}^2y_{14}
 {+}y_7^2y_8^2y_{11}^2y_{13}
 {+}y_7^3y_8y_{11}y_{12}y_{13}
 {+}y_7^3y_8y_{11}^2y_{14}
  \\ & \qquad
 {+}y_7^4y_{12}^2y_{13}
 {+}y_7^4y_{11}y_{12}y_{14}
 {+}y_7^4y_{11}^2y_{15}
 {+}y_7^4y_8^3y_{13}
 {+}y_7^5y_{15}^2
 {+}y_7^5y_8^2y_{14}
  \\ & \qquad
 {+}y_7^6y_8y_{15}
 {+}y_6y_{10}y_{11}^2y_{13}y_{14}
 {+}y_6y_{10}^2y_{13}^3
 {+}y_6y_{10}^2y_{11}y_{13}y_{15}
  \\ & \qquad
 {+}y_6y_7y_{11}y_{13}^2y_{15}
 {+}y_6y_7^2y_8y_{11}y_{13}^2
 {+}y_6y_7^2y_8y_{11}^2y_{15}
 {+}y_6y_7^3y_{11}y_{13}y_{14}
  \\ & \qquad
 {+}y_6y_7^3y_{11}y_{12}y_{15}
 {+}y_6y_7^4y_8^2y_{15}
 {+}y_6^2y_8^2y_{11}y_{13}^2
 {+}y_6^2y_7^2y_{13}^3
 {+}y_6^2y_7^2y_{11}y_{14}^2
  \\ & \qquad
 {+}y_6^2y_7^2y_{11}y_{13}y_{15}
 {+}y_6^3y_{11}^2y_{12}y_{13}
 {+}y_6^3y_{10}y_{11}y_{13}^2
 {+}y_6^3y_{10}y_{11}^2y_{15}
  \\ & \qquad
 {+}y_6^3y_{10}^2y_{13}y_{14}
 {+}y_6^3y_8y_{13}^3
 {+}y_6^3y_7y_8^2y_{11}y_{13}
 {+}y_6^3y_7^2y_8y_{12}y_{13}
 {+}y_6^4y_{13}y_{14}^2
  \\ & \qquad
 {+}y_6^4y_8^2y_{12}y_{13}
 {+}y_6^4y_8^2y_{11}y_{14}
 {+}y_6^4y_8^2y_{10}y_{15}
 {+}y_6^4y_7y_8y_{13}^2
 {+}y_6^4y_7y_8y_{11}y_{15}
  \\ & \qquad
 {+}y_6^5y_8y_{13}y_{14}
 {+}y_4y_{11}^2y_{13}^3
 {+}y_4y_{10}y_{11}y_{13}^2y_{14}
 {+}y_4y_7^2y_{11}^2y_{12}y_{13}
  \\ & \qquad
 {+}y_4y_7^2y_8y_{11}y_{13}y_{15}
 {+}y_4y_7^3y_{13}^2y_{14}
 {+}y_4y_7^3y_{11}y_{14}y_{15}
 {+}y_4y_7^3y_8^2y_{11}y_{13}
  \\ & \qquad
 {+}y_4y_7^4y_8y_{11}y_{14}
 {+}y_4y_7^5y_{12}y_{14}
 {+}y_4y_7^5y_{11}y_{15}
 {+}y_4y_6y_7y_{11}^2y_{13}^2
  \\ & \qquad
 {+}y_4y_6y_7^2y_{13}^2y_{15}
 {+}y_4y_6y_7^2y_{11}y_{15}^2
 {+}y_4y_6y_7^3y_8y_{11}y_{15}
 {+}y_4y_6y_7^4y_{13}y_{14}
  \\ & \qquad
 {+}y_4y_6y_7^4y_{12}y_{15}
 {+}y_4y_6^2y_{11}^2y_{13}y_{14}
 {+}y_4y_6^2y_{10}y_{11}y_{13}y_{15}
 {+}y_4y_6^2y_7^2y_8y_{13}y_{14}
  \\ & \qquad
 {+}y_4y_6^2y_7^3y_{14}^2
 {+}y_4y_6^2y_7^3y_{13}y_{15}
 {+}y_4y_6^3y_7y_{11}y_{12}y_{13}
 {+}y_4y_6^4y_{11}y_{13}^2
  \\ & \qquad
 {+}y_4y_6^4y_{11}^2y_{15}
 {+}y_4y_6^4y_{10}y_{13}y_{14}
 {+}y_4y_6^4y_7y_{15}^2
 {+}y_4^2y_{10}y_{11}^2y_{12}y_{13}
  \\ & \qquad
 {+}y_4^2y_{10}^2y_{11}y_{13}^2
 {+}y_4^2y_{10}^2y_{11}^2y_{15}
 {+}y_4^2y_{10}^3y_{13}y_{14}
 {+}y_4^2y_8y_{11}y_{12}y_{13}^2
  \\ & \qquad
 {+}y_4^2y_8y_{11}^2y_{13}y_{14}
 {+}y_4^2y_8y_{10}y_{11}y_{13}y_{15}
 {+}y_4^2y_7y_{11}y_{13}^3
 {+}y_4^2y_7y_{11}y_{12}y_{13}y_{14}
  \\ & \qquad
 {+}y_4^2y_7y_{11}^2y_{14}^2
 {+}y_4^2y_7y_{11}^2y_{13}y_{15}
 {+}y_4^2y_7^3y_{11}y_{12}y_{13}
 {+}y_4^2y_6y_{11}y_{13}^2y_{14}
  \\ & \qquad
 {+}y_4^2y_6y_{11}y_{12}y_{13}y_{15}
 {+}y_4^2y_6y_{10}y_{13}^2y_{15}
 {+}y_4^2y_6y_{10}y_{11}y_{15}^2
  \\ & \qquad
 {+}y_4^2y_6y_7y_8y_{11}y_{12}y_{13}
 {+}y_4^2y_6y_7^2y_{12}^2y_{13}
 {+}y_4^2y_6y_7^2y_{11}y_{13}^2
 {+}y_4^2y_6^2y_{10}y_{11}^2y_{13}
  \\ & \qquad
 {+}y_4^2y_6^2y_{10}^2y_{12}y_{13}
 {+}y_4^2y_6^2y_{10}^2y_{11}y_{14}
 {+}y_4^2y_6^2y_{10}^3y_{15}
 {+}y_4^2y_6^2y_8y_{11}^2y_{15}
  \\ & \qquad
 {+}y_4^2y_6^2y_8y_{10}y_{13}y_{14}
 {+}y_4^2y_6^2y_7y_{11}y_{12}y_{15}
 {+}y_4^2y_6^3y_{12}y_{13}y_{14}
 {+}y_4^2y_6^3y_{13}^3
  \\ & \qquad
 {+}y_4^3y_{11}y_{13}y_{14}y_{15}
 {+}y_4^3y_6y_{11}^2y_{12}y_{13}
 {+}y_4^3y_6y_{10}y_{11}y_{13}^2
 {+}y_4^3y_6y_{10}y_{11}^2y_{15}
  \\ & \qquad
 {+}y_4^3y_6y_{10}^2y_{13}y_{14}
 {+}y_4^3y_6y_8y_{11}y_{13}y_{15}
 {+}y_4^3y_6y_7y_{12}y_{13}y_{15}
 {+}y_4^3y_6^2y_{13}y_{14}^2
  \\ & \qquad
 {+}y_4^3y_6^2y_{13}^2y_{15}
 {+}y_4^3y_7y_8y_{11}y_{13}y_{14}
 {+}y_4^3y_7^2y_{11}y_{14}^2
 {+}y_4^3y_7^2y_{11}y_{13}y_{15}
  \\ & \qquad
 {+}y_4^4y_{10}y_{12}y_{13}y_{14}
 {+}y_4^4y_{10}y_{12}^2y_{15}
 {+}y_4^4y_{11}y_{12}y_{13}^2
 {+}y_4^4y_{11}y_{12}^2y_{14}
  \\ & \qquad
 {+}y_4^4y_{11}^2y_{12}y_{15}
 {+}y_4^4y_{12}^3y_{13}
 {+}y_4^4y_8y_{13}y_{14}^2
 {+}y_4^4y_8y_{13}^2y_{15}
 {+}y_4^4y_8y_{11}y_{15}^2
  \\ & \qquad
 {+}y_4^4y_7y_{14}^3
 {+}y_4^4y_7y_{12}y_{15}^2
 {+}y_4^4y_6y_{14}^2y_{15}
 {+}y_4^4y_6y_{13}y_{15}^2
 {+}y_4^5y_{15}^3
  ,
\end{align*}
where the first term of
$R_i$ is the leading term.

In order to calculate the Stiefel--Whitney classes of the spin
representation concretely, we need some facts about the spin
representations.
\begin{notation}
  $f_n \co \Spin(n) \longrightarrow \Spin(n+1)$ is the natural inclusion
  map.
\end{notation}
As is well known (see for example Adams \cite[Proposition 4.4]{adams2}),
we have that
\begin{align*}
  & f_{2k-1}^* (\Delta_{2k}^{\pm})_{\mathbb C} =
  (\Delta_{2k-1})_{\mathbb C}, \\
  & f_{2k}^* (\Delta_{2k+1})_{\mathbb C} = (\Delta_{2k}^+)_{\mathbb C}
  \oplus (\Delta_{2k}^-)_{\mathbb C}.
\end{align*}
Then, using these and \fullref{rem:iso},
we have the following:
\begin{equation}
  \label{eq:spin_rep}
  \begin{split}
    & f_{8m}^* \Delta_{8m+1} = \Delta_{8m}^+ \oplus \Delta_{8m}^-, \\
    & f_{8m+1}^* \Delta_{8m+2}^{\pm} = (\Delta_{8m+1})_{\mathbb C}, \\
    & f_{8m+2}^* (\Delta_{8m+3})_{\mathbb C} = \Delta_{8m+2}^+ \oplus
    \Delta_{8m+2}^-, \\
    & f_{8m+3}^* \Delta_{8m+4}^{\pm} = \Delta_{8m+3}, \\
    & f_{8m+4}^* \Delta_{8m+5} = \Delta_{8m+4}^+ \oplus \Delta_{8m+4}^-,
    \\
    & f_{8m+5}^* \Delta_{8m+6}^{\pm} = (\Delta_{8m+5})_{\mathbb C},
    \\
    & f_{8m+6}^* \Delta_{8m+7} = (\Delta_{8m+6}^+)_{\mathbb R}
    = (\Delta_{8m+6}^-)_{\mathbb R}, \\
    & f_{8m+7}^* \Delta_{8m+8}^{\pm} = \Delta_{8m+7}. 
  \end{split}
\end{equation}

\begin{notation}
  $\lambda^i \co G(n) \to G(\tbinom{n}{i})$ is the $i$--th exterior
  representation, where $G = O, SO, U, SU$.
\end{notation}

For the usual inclusion map
\begin{displaymath}
  r_n \co SU(n) \to SO(2n),
\end{displaymath}
there is a covering map
\begin{displaymath}
  \tilde{r}_n \co SU(n) \to \Spin(2n),
\end{displaymath}
since $SU(n)$ is simply connected.
According to Atiyah, Bott and Shapiro \cite{ABS},
there are isomorphisms of the representations
\begin{align}
  \label{eq:plus}
  & \tilde{r}_n^* (\Delta_{2n}^+)_{\mathbb C}
  = \sum_{i=0}^{[n/2]} \lambda^{2i}, \\
  \label{eq:minus}
  & \tilde{r}_n^* (\Delta_{2n}^-)_{\mathbb C}
  = \sum_{i=0}^{[(n-1)/2]} \lambda^{2i+1}.
\end{align}

\section{The Stiefel--Whitney classes of the second exterior
  representation}
\label{sec:2nd_exterior_rep}

In this section, we calculate the Stiefel--Whitney classes of the
representation
\begin{displaymath}
  \lambda_{15}^2 \co \Spin(15) \longrightarrow SO(15) \longrightarrow
  O(15) \longrightarrow O(105)
\end{displaymath}
induced from the second exterior representation $\lambda^2 \co O(15) \to
O(105)$.

Let $\nu$ be the nontrivial real representation of one dimension
\begin{displaymath}
  \nu \co \mathbb Z/2 \overset{\cong}{\longrightarrow} O(1),
\end{displaymath}
and $\nu_i$ a composition map
\begin{displaymath}
  \nu_i \co H_n \cong (\mathbb Z/2)^n \overset{p_i}{\longrightarrow}
  \mathbb Z/2 \overset{\nu}{\longrightarrow} O(1),
\end{displaymath}
where $p_i$ is the $i$--th projection map.
Then it is easy to show that the induced representation
\begin{displaymath}
  i_n^* \lambda^2 \co H_n \overset{i_n}{\longrightarrow} O(n)
  \overset{\lambda^2}{\longrightarrow}
  O\big(\tbinom{n}{2}\big)
\end{displaymath}
is isomorphic to ${\sum_{1\leq j<k\leq n} \nu_j\nu_k}$.
Using \fullref{prop:BH}, we obtain
\begin{displaymath}
  Bi_n^* w(\lambda^2) = \prod_{1 \leq j < k \leq n} (1 + t_j + t_k),
\end{displaymath}
since the total Stiefel--Whitney class of the representation
$\nu_j\nu_k$ is given by $1 + t_j + t_k$.
Then we can write
\begin{displaymath}
  \varphi_n^2(w_1, w_2, \ldots, w_n)
  = \prod_{1 \leq j < k \leq n} (1 + t_j + t_k),
\end{displaymath}
where $w_i$ is the $i$--th elementary symmetric polynomial.
Since it is quite difficult to calculate $\varphi_{15}^2$ directly
by using GAP because of the capacity of the computer, we need to
improve the algorithm to expand symmetric polynomials in terms of
the elementary symmetric polynomials, that is, to calculate
$\varphi_n^2$ by induction on $n$.
Obviously $\varphi_2^2(w_1, w_2) = 1 + w_1$.
Suppose that $\varphi_{n-1}^2$ is obtained.
Let $\bar{w}_i$ be the $i$--th elementary symmetric polynomial of
$t_j$ $(1 \leq j \leq n-1)$.
Then we have
\begin{displaymath}
  w_i =
  \begin{cases}
    \bar{w}_1 + t_n, & \textrm{ if } i = 1, \\
    \bar{w}_i + \bar{w}_{i-1}t_n & \textrm{ if } 1 < i < n, \\
    \bar{w}_{n-1}t_n & \textrm{ if } i = n.
  \end{cases}
\end{displaymath}
We define a function $\psi_n^2$ by
\begin{displaymath}
  \psi_n^2(\bar{w}_1, \ldots, \bar{w}_{n-1}, t_n)
  = \prod_{1 \leq j < k \leq n-1} (1 + t_j + t_k)
  \prod_{i=1}^{n-1} (1 + t_i + t_n),
\end{displaymath}
which is equal to $\varphi_n^2(w_1, \ldots, w_n)$ as a polynomial of $t_i$.
Then we obtain
\begin{displaymath}
  \psi_n^2(\bar{w}_1, \ldots, \bar{w}_{n-1}, t_n)
  = \varphi_{n-1}^2(\bar{w}_1, \bar{w}_2, \ldots, \bar{w}_{n-1})
  \sum_{k=0}^{n-1} \Big( t_n^k \sum_{l=0}^{n-1-k} \binom{n-1-l}{k}
    \bar{w}_l \Big).
\end{displaymath}
We define $\varphi_{n,i}^2(w_1, \ldots, w_{n-1})$ for $i \geq 0$ by
the equation
\begin{displaymath}
  \varphi_n^2(w_1, \ldots, w_n)
  = \sum_m \varphi_{n,m}^2(w_1, \ldots, w_{n-1}) w_n^m,
\end{displaymath}
where it is easy to see that there holds the following identity
\begin{displaymath}
  \varphi_{n,0}^2(\xi_1, \ldots, \xi_{n-1})
  = \psi_n^2(\xi_1, \ldots, \xi_{n-1}, 0),
\end{displaymath}
as polynomials for any invariant element $\xi_i$.
Now we calculate $\varphi_{n,m}^2$ by induction on $m$.
Assume that $\varphi_{n,l}^2$ is obtained for $1 \leq l \leq m-1$.
Put
\begin{align*}
  \psi_{n,l}^2(\bar{w}_1, \ldots, \bar{w}_{n-1}, t_n)
  &= \varphi_{n,l}^2(w_1, \ldots, w_{n-1})\\
\tag*{\hbox{and}}
  \chi_{n,m}^2(\bar{w}_1,
  .\, .\, ,
  \bar{w}_{n-1}, t_n)
  & = \{ \psi_n^2(\bar{w}_1,
  .\, .\, ,
  \bar{w}_{n-1}, t_n)
  \\ & \qquad
  - \sum\limits_{l=0}^{m-1} \psi_{n,l}^2(\bar{w}_1,
  .\, .\, ,
  \bar{w}_{n-1}, t_n)
  \bar{w}_{n-1}^lt_n^l \}
  /\bar{w}_{n-1}^mt_n^m.\\
\tag*{\hbox{Then}}
  \varphi_{n,m}^2(\xi_1, \ldots, \xi_{n-1})
  &= \chi_{n,m}^2(\xi_1, \ldots, \xi_{n-1}, 0),
\end{align*}
which gives $\varphi_n^2(w_1, \ldots, w_n)$ by the above equality.
By induction on $n$, namely, using the improved algorithm,
we obtain the following theorem.
\begin{thm}
  The Stiefel--Whitney classes of degree $2^i$ for $0 \leq i \leq 6$
  of the induced
  representation $\lambda_{15}^2 \co \Spin(15) \to O(105)$ are given as
  follows:
  \begin{align*}
    w_1(\lambda_{15}^2) & = 0, \\
    w_2(\lambda_{15}^2) & = 0, \\
    w_4(\lambda_{15}^2) & =
    y_4
    , \\
    w_8(\lambda_{15}^2) & =
    y_8
   {+}y_4^2
    , \\
    w_{16}(\lambda_{15}^2) & =
    y_8^2
   {+}y_4^2y_8
   {+}y_4y_6^2
    , \\
    w_{32}(\lambda_{15}^2) & =
    y_4y_{14}^2
   {+}y_6y_{13}^2
   {+}y_8y_{12}^2
   {+}y_{10}^2y_{12}
   {+}y_{10}y_{11}^2
   {+}y_4^2y_{11}y_{13}
   {+}y_6^2y_7y_{13}
    \\ & \ \ 
   {+}y_6^2y_{10}^2
   {+}y_7^3y_{11}
   {+}y_8^4
   {+}y_4^8
    , \\
    w_{64}(\lambda_{15}^2) & =
    y_4y_{15}^4
   {+}y_8y_{14}^4
   {+}y_4^2y_{14}^4
   {+}y_4y_6y_{13}^3y_{15}
    \\ & \ \ 
    {+}y_4y_7y_{12}y_{13}^2y_{15}
   {+}y_4y_8y_{11}y_{13}^2y_{15}
    \\ & \ \ 
   {+}y_4y_8y_{13}^4
   {+}y_4y_{10}^2y_{12}y_{13}y_{15}
   {+}y_4y_{10}^2y_{13}^2y_{14}
   {+}y_4y_{10}y_{11}^2y_{13}y_{15}
    \\ & \ \ 
   {+}y_4y_{10}y_{11}y_{13}^3
   {+}y_4y_{11}^3y_{12}y_{15}
   {+}y_4y_{11}^2y_{12}y_{13}^2
   {+}y_6^2y_{13}^4
   {+}y_7y_{11}^4y_{13}
    \\ & \ \ 
   {+}y_8^2y_{12}^4
   {+}y_8y_{10}^3y_{11}y_{15}
   {+}y_8y_{10}^2y_{11}^2y_{14}
   {+}y_8y_{10}y_{11}^3y_{13}
   {+}y_8y_{11}^4y_{12}
    \\ & \ \ 
   {+}y_{10}^4y_{11}y_{13}
   {+}y_{10}^2y_{11}^4
   {+}y_4^3y_{10}y_{12}y_{15}^2
   {+}y_4^3y_{10}y_{13}y_{14}y_{15}
   {+}y_4^3y_{11}^2y_{15}^2
    \\ & \ \ 
   {+}y_4^3y_{11}y_{12}y_{14}y_{15}
   {+}y_4^3y_{12}^2y_{13}y_{15}
    \\ & \ \ 
   {+}y_4^3y_{12}y_{13}^2y_{14}
   {+}y_4^2y_6y_{10}y_{12}y_{13}y_{15}
    \\ & \ \ 
   {+}y_4^2y_6y_{10}y_{13}^2y_{14}
   {+}y_4^2y_6y_{11}^2y_{13}y_{15}
   {+}y_4^2y_6y_{11}y_{13}^3
   {+}y_4^2y_7^2y_{12}y_{15}^2
    \\ & \ \ 
   {+}y_4^2y_7y_8y_{13}^2y_{15}
   {+}y_4^2y_7y_{11}^2y_{12}y_{15}
   {+}y_4^2y_7y_{11}y_{12}y_{13}^2
   {+}y_4^2y_8^2y_{13}^2y_{14}
    \\ & \ \ 
   {+}y_4^2y_8y_{10}^2y_{14}^2
   {+}y_4^2y_8y_{10}y_{11}y_{12}y_{15}
   {+}y_4^2y_8y_{11}^2y_{12}y_{14}
   {+}y_4^2y_8y_{12}^4
    \\ & \ \ 
   {+}y_4^2y_{10}^3y_{11}y_{15}
   {+}y_4^2y_{10}^3y_{13}^2
   {+}y_4^2y_{10}^2y_{11}^2y_{14}
       \\ & \ \ 
   {+}y_4^2y_{10}y_{11}^3y_{13}
   {+}y_4y_6^3y_{12}y_{15}^2
    \\ & \ \ 
   {+}y_4y_6^3y_{13}y_{14}y_{15}
   {+}y_4y_6^2y_7y_{12}y_{14}y_{15}
   {+}y_4y_6^2y_8y_{10}y_{15}^2
    \\ & \ \ 
   {+}y_4y_6^2y_8y_{11}y_{14}y_{15}
   {+}y_4y_6^2y_8y_{12}y_{13}y_{15}
   {+}y_4y_6^2y_8y_{13}^2y_{14}
    \\ & \ \ 
   {+}y_4y_6^2y_{10}^2y_{13}y_{15}
   {+}y_4y_6^2y_{10}^2y_{14}^2
   {+}y_4y_6^2y_{11}^3y_{15}
       \\ & \ \ 
   {+}y_4y_6^2y_{12}^4
   {+}y_4y_6y_7^2y_{10}y_{15}^2
    \\ & \ \ 
   {+}y_4y_6y_7y_8y_{13}^3
   {+}y_4y_6y_8^2y_{10}y_{13}y_{15}
   {+}y_4y_6y_8y_{10}^2y_{11}y_{15}
    \\ & \ \ 
   {+}y_4y_6y_8y_{10}y_{11}^2y_{14}
   {+}y_4y_6y_8y_{11}^3y_{13}
   {+}y_4y_7^3y_{12}^2y_{15}
   {+}y_4y_7^2y_8y_{12}y_{13}^2
    \\ & \ \ 
   {+}y_4y_7^2y_{10}^2y_{13}^2
   {+}y_4y_7^2y_{11}^2y_{12}^2
   {+}y_4y_7y_8^2y_{11}^2y_{15}
   {+}y_4y_7y_8^2y_{11}y_{13}^2
    \\ & \ \ 
   {+}y_4y_7y_8y_{11}^3y_{12}
   {+}y_4y_8^4y_{14}^2
   {+}y_4y_8^2y_{11}^4
   {+}y_4y_{10}^6
   {+}y_6^2y_7^2y_8y_{15}^2
    \\ & \ \ 
   {+}y_6^2y_7^2y_{12}y_{13}^2
   {+}y_6^2y_8y_{11}^4
   {+}y_6^2y_{10}^4y_{12}
\\& \ \ 
   {+}y_6^2y_{10}^3y_{11}^2
   {+}y_6y_7^2y_{11}^4
   {+}y_6y_8^4y_{13}^2
    \\ & \ \ 
   {+}y_7^4y_8y_{14}^2
   {+}y_7^4y_{10}y_{13}^2
   {+}y_8^5y_{12}^2
   {+}y_8^4y_{10}^2y_{12}
   {+}y_8^4y_{10}y_{11}^2
   {+}y_8^3y_{10}^4
    \\ & \ \ 
   {+}y_4^4y_7y_{11}y_{15}^2
   {+}y_4^4y_7y_{13}^2y_{15}
   {+}y_4^4y_7y_{13}y_{14}^2
  \\ & \ \
   {+}y_4^4y_{11}y_{12}^2y_{13}
   {+}y_4^3y_6y_8^2y_{15}^2
    \\ & \ \ 
   {+}y_4^3y_6y_{11}^3y_{13}
   {+}y_4^3y_7y_8^2y_{14}y_{15}
   {+}y_4^3y_7y_8y_{11}^2y_{15}
   {+}y_4^3y_7y_{10}^3y_{15}
    \\ & \ \ 
   {+}y_4^3y_7y_{10}^2y_{11}y_{14}
   {+}y_4^3y_7y_{11}^3y_{12}
   {+}y_4^2y_6^2y_7^2y_{15}^2
   {+}y_4^2y_6^2y_8^2y_{14}^2
    \\ & \ \ 
   {+}y_4^2y_6^2y_{10}^2y_{11}y_{13}
   {+}y_4^2y_6^2y_{10}^2y_{12}^2
\\ & \ \ 
   {+}y_4^2y_6y_7y_8^2y_{12}y_{15}
   {+}y_4^2y_6y_7y_8^2y_{13}y_{14}
    \\ & \ \ 
   {+}y_4^2y_6y_7y_8y_{10}^2y_{15}
   {+}y_4^2y_6y_7y_8y_{11}^2y_{13}
   {+}y_4^2y_7^4y_{13}y_{15}
   {+}y_4^2y_7^4y_{14}^2
    \\ & \ \ 
   {+}y_4^2y_7^3y_{11}^2y_{13}
   {+}y_4^2y_7^2y_8^2y_{11}y_{15}
   {+}y_4^2y_7^2y_8^2y_{12}y_{14}
   {+}y_4^2y_7^2y_8^2y_{13}^2
    \\ & \ \ 
   {+}y_4^2y_7^2y_8y_{10}^2y_{14}
   {+}y_4^2y_7^2y_8y_{11}^2y_{12}
   {+}y_4^2y_8^4y_{11}y_{13}
   {+}y_4y_6^5y_{15}^2
    \\ & \ \ 
   {+}y_4y_6^4y_7y_{14}y_{15}
   {+}y_4y_6^4y_{10}y_{13}^2
   {+}y_4y_6^2y_7y_8^2y_{10}y_{15}
   {+}y_4y_6^2y_7y_8^2y_{11}y_{14}
    \\ & \ \ 
   {+}y_4y_6y_7^4y_{11}y_{15}
   {+}y_4y_6y_7^2y_8^2y_{11}y_{13}
   {+}y_4y_7^5y_{11}y_{14}
   {+}y_4y_7^3y_8^3y_{15}
    \\ & \ \ 
   {+}y_4y_7^3y_8^2y_{11}y_{12}
   {+}y_6^6y_{14}^2
   {+}y_6^5y_7y_{12}y_{15}
   {+}y_6^5y_7y_{13}y_{14}
   {+}y_6^5y_{10}y_{11}y_{13}
    \\ & \ \ 
   {+}y_6^4y_7^2y_{11}y_{15}
   {+}y_6^4y_7^2y_{12}y_{14}
   {+}y_6^4y_7^2y_{13}^2
   {+}y_6^4y_7y_{11}^3
   {+}y_6^4y_8^2y_{11}y_{13}
    \\ & \ \ 
   {+}y_6^4y_8^2y_{12}^2
   {+}y_6^3y_7y_8^3y_{15}
   {+}y_6^2y_7^2y_8^3y_{14}
   {+}y_6^2y_7y_8^4y_{13}
   {+}y_6^2y_8^4y_{10}^2
    \\ & \ \ 
   {+}y_6y_7^5y_8y_{15}
   {+}y_6y_7^3y_8^3y_{13}
   {+}y_7^7y_{15}
   {+}y_7^6y_8y_{14}
   {+}y_7^5y_8^2y_{13}
   {+}y_7^4y_8^3y_{12}
    \\ & \ \ 
   {+}y_7^3y_8^4y_{11}
   {+}y_4^6y_{13}^2y_{14}
   {+}y_4^5y_7y_{11}^2y_{15}
   {+}y_4^5y_7y_{11}y_{13}^2
   {+}y_4^4y_6^2y_{11}^2y_{14}
    \\ & \ \ 
   {+}y_4^4y_6y_7y_{10}^2y_{15}
   {+}y_4^4y_6y_7y_{11}^2y_{13}
   {+}y_4^4y_7^2y_8y_{13}^2
   {+}y_4^4y_7^2y_{11}^2y_{12}
    \\ & \ \ 
   {+}y_4^3y_6^2y_7^2y_{13}^2
   {+}y_4^3y_7^4y_{11}y_{13}
   {+}y_4^2y_6^5y_{13}^2
   \\ & \ \ 
   {+}y_4^2y_6^4y_7y_{10}y_{15}
   {+}y_4^2y_6^4y_7y_{11}y_{14}
    \\ & \ \ 
   {+}y_4^2y_6^4y_7y_{12}y_{13}
   {+}y_4^2y_6^4y_8y_{11}y_{13}
   {+}y_4y_6^6y_{11}y_{13}
   {+}y_4y_6^4y_7^3y_{15}
    \\ & \ \ 
   {+}y_4y_6^2y_7^5y_{13}
   {+}y_4y_7^7y_{11}
   {+}y_4y_7^4y_8^4
   {+}y_6^7y_7y_{15}
   {+}y_6^6y_7^2y_{14}
   {+}y_6^6y_8y_{10}^2
    \\ & \ \ 
   {+}y_6^5y_7^3y_{13}
   {+}y_6^4y_7^4y_{12}
   {+}y_6^4y_8^5
   {+}y_7^8y_8
   {+}y_4^7y_6y_{15}^2
   {+}y_4^7y_7y_{14}y_{15}
    \\ & \ \ 
   {+}y_4^6y_6^2y_{14}^2
   {+}y_4^6y_6y_7y_{12}y_{15}
   {+}y_4^6y_6y_7y_{13}y_{14}
\\ & \ \ 
   {+}y_4^6y_7^2y_{11}y_{15}
   {+}y_4^6y_7^2y_{12}y_{14}
    \\ & \ \ 
   {+}y_4^6y_{10}^4
   {+}y_4^5y_6^2y_7y_{10}y_{15}
   {+}y_4^5y_6^2y_7y_{11}y_{14}
\\ & \ \ 
   {+}y_4^5y_6y_7^2y_{11}y_{13}
   {+}y_4^5y_7^3y_8y_{15}
    \\ & \ \ 
   {+}y_4^5y_7^3y_{11}y_{12}
   {+}y_4^4y_6^4y_{12}^2
   {+}y_4^4y_6^3y_7y_8y_{15}
\\ & \ \ 
   {+}y_4^4y_6^2y_7^2y_8y_{14}
   {+}y_4^4y_6y_7^3y_8y_{13}
    \\ & \ \ 
   {+}y_4^4y_7^5y_{13}
   {+}y_4^4y_7^4y_8y_{12}
   {+}y_4^4y_8^6
   {+}y_4^2y_6^6y_7y_{13}
   {+}y_4^2y_6^6y_{10}^2
   {+}y_4^2y_6^4y_7^3y_{11}
    \\ & \ \ 
   {+}y_4^2y_6^4y_8^4
   {+}y_4^2y_7^8
   {+}y_6^8y_8^2
   {+}y_4^9y_{14}^2
   {+}y_4^8y_7y_{10}y_{15}
   {+}y_4^8y_7y_{11}y_{14}
    \\ & \ \ 
   {+}y_4^8y_7y_{12}y_{13}
   {+}y_4^8y_8y_{11}y_{13}
   {+}y_4^8y_8y_{12}^2
   {+}y_4^8y_{10}^2y_{12}
   {+}y_4^8y_{10}y_{11}^2
    \\ & \ \ 
   {+}y_4^7y_7^3y_{15}
   {+}y_4^7y_7^2y_{11}^2
   {+}y_4^6y_6^3y_7y_{15}
   {+}y_4^6y_6^2y_7^2y_{14}
   {+}y_4^6y_6^2y_8y_{10}^2
    \\ & \ \ 
   {+}y_4^6y_6y_7^3y_{13}
   {+}y_4^6y_7^4y_{12}
   {+}y_4^6y_8^5
   {+}y_4^5y_6^4y_7y_{13}
   {+}y_4^5y_6^4y_{10}^2
   {+}y_4^5y_6^2y_7^3y_{11}
    \\ & \ \ 
   {+}y_4^5y_6^2y_8^4
   {+}y_4^2y_6^8y_8
   {+}y_4y_6^{10}
    .
  \end{align*}
\end{thm}

\begin{rem}
  The total Stiefel--Whitney classes of the induced representation
  $$\lambda_n^2 \co \Spin(n) \to O\big(\tbinom{n}{2}\big)$$ for $3 \leq n \leq 9$
  are given as follows:
  \begin{align*}
    w(\lambda_3^2) & = 1, \\
    w(\lambda_4^2) & = 1, \\
    w(\lambda_5^2) & = 1 + y_4, \\
    w(\lambda_6^2) & = (1 + y_4 + y_6)^2, \\
    w(\lambda_7^2) & = (1 + y_4 + y_6 + y_7)^3, \\
    w(\lambda_8^2) & = (1 + y_4 + y_6 + y_7)^4, \\
    w(\lambda_9^2) & = (1 + y_4 + y_6 + y_7 + y_8)(1 + y_4 + y_6 + y_7)^4.
  \end{align*}
\end{rem}

\section{The Stiefel--Whitney classes of the spin representation}
\label{sec:spin_rep}

In this section, by making use of \fullref{thm:quillen}, we calculate
the Stiefel--Whitney classes of the spin representations
$(\Delta_n)_{\mathbb R}$ and $(\Delta_n^{\pm})_{\mathbb R}$ for $n \leq 15$.
Firstly, one can obtain them for
$n \leq 7$ as follows:
\begin{align*}
  & w((\Delta_3)_{\mathbb R}) = p(\Delta_3) = 1 + u_4, \\
  & w((\Delta_4^+)_{\mathbb R}) = p(\Delta_4^+) = 1 + u_4, \\
  & w((\Delta_4^-)_{\mathbb R}) = p(\Delta_4^-) = 1 + y_4 + u_4, \\
  & w((\Delta_5)_{\mathbb R}) = p(\Delta_5) = 1 + y_4 + u_8, \\
  & w((\Delta_6)_{\mathbb R}) = c(\Delta_6) = 1 + y_4 + y_6 + u_8, \\
  & w(\Delta_7) = 1 + y_4 + y_6 + y_7 + u_8.
\end{align*}
Secondly, it also follows from \fullref{thm:quillen} that the total
Stiefel--Whitney class of the spin representation $\Delta_8^+$ is given
by
\begin{displaymath}
  w(\Delta_8^+) = 1 + y_4 + y_6 + y_7 + u_8.
\end{displaymath}
According to Adams \cite{adams2}, the outer automorphism group
$\mathrm{Out}(\Spin(8))$ of $\Spin(8)$ is isomorphic to the symmetric
group $\Sigma_3$ of degree 3 which  acts on the set of the
representations $\lambda^1$, $\Delta_8^+$ and $\Delta_8^-$, where
$\lambda^1$ is the natural projection $\Spin(8) \to SO(8)$.
Then there is an automorphism $\sigma \co \Spin(8) \to \Spin(8)$ such
that $\sigma^*(\lambda^1) = \smash{\Delta_8^-}$, $\sigma^*(\smash{\Delta_8^+}) =
\lambda^1$ and $\sigma^*(\smash{\Delta_8^-}) = \smash{\Delta_8^+}$.
Since $Bf_7^*w(\smash{\Delta_8^-}) = w(\Delta_7)$,
we can write as follows:
\begin{displaymath}
  w(\Delta_8^-) = 1 + y_4 + y_6 + y_7 + u_8 + a_1y_8,
\end{displaymath}
where $a_1 \in \mathbb Z/2$. Then we have
\begin{align*}
  B\sigma^* y_8 & = B\sigma^* w_8(\lambda_8^1) = w_8 (\Delta_8^-)
  = u_8 + a_1y_8, \\
  B\sigma^* u_8 & = B\sigma^* w_8(\Delta_8^+) = w_8 (\lambda_8^1) = y_8,\\
\tag*{\hbox{and}}
  u_8 &= w_8(\Delta_8^+)
  = B\sigma^* w_8(\Delta_8^-)
  = y_8 + a_1(u_8 + a_1y_8).
\end{align*}
Thus we obtain $a_1 = 1$, and the total Stiefel--Whitney class is
given by
\begin{displaymath}
  w(\Delta_8^-) = 1 + y_4 + y_6 + y_7 + (u_8 + y_8).
\end{displaymath}
Since the induced representation $f_8^*\Delta_9$ is
isomorphic to $\Delta_8^+ \oplus \Delta_8^-$ by \eqref{eq:spin_rep},
the total Stiefel--Whitney class of $f_8^*\Delta_9$ is given by
\begin{align*}
  w(f_8^*\Delta_9) & = w(\Delta_8^+)w(\Delta_8^-) \\
  & = 1 + (y_8 + y_4^2) + (y_4y_8 + y_6^2) + (y_6y_8 + y_7^2) + y_7y_8
  + (u_8^2 + y_8u_8).
\end{align*}
Since $Bf_8^* \co H^*(B\Spin(9);\mathbb Z/2) \to
H^*(B\Spin(8);\mathbb Z/2)$ is a monomorphism,
the total Stiefel--Whitney class of the spin representation
$\Delta_9$ is given by
\begin{equation}
  \label{eq:delta9}
  w(\Delta_9) = 1 + (y_8 + y_4^2) + (y_4y_8 + y_6^2)
  + (y_6y_8 + y_7^2)
  + y_7y_8 + u_{16}.
\end{equation}
Recall from \eqref{eq:plus} that we have $\tilde{r}_5^*\Delta_{10}^+
= 1 \oplus \lambda^2 \oplus \lambda^4$,
where $\lambda^i$ is the $i$--th exterior representation $SU(5) \to
SU(\smash{\tbinom{5}{i}})$.
In a similar way to the Stiefel--Whitney classes of $i$--the exterior representations,
we can calculate the mod 2 Chern classes of the $i$--th exterior
representations $\lambda^i$ using the Borel--Hirzebruch method;
the mod 2 total Chern classes are given by
\begin{align*}
  c(\lambda^2) & = 1 + c_2 + c_3 + (c_4 + c_2^2) + c_5 + (c_3^2 + c_2^3)
  + c_2^2c_3 + (c_3c_5 + c_2^2c_4 + c_2c_3^2)
  \\ & \quad
  + (c_3^3 + c_2^2c_5)
  + (c_5^2 + c_2c_3c_5 + c_3^2c_4), \\
  c(\lambda^4) & = 1 + c_2 + c_3 + c_4 + c_5,
\end{align*}
where $c_i \in H^*(BSU(5);\mathbb Z/2)$ is the mod 2 $i$--th Chern class.
By \eqref{eq:whitney_sum}, we obtain
\begin{displaymath}
  w_{16}(\tilde{r}_5^*(\Delta_{10}^+)_{\mathbb R})
  = c_8(\tilde{r}_5^*\Delta_{10}^+)
  = c_3c_5 + c_4^2 + c_2^4.
\end{displaymath}
Since $B\tilde{r}_5^*\co H^{16}(B\Spin(10);\mathbb Z/2) \to
H^{16}(BSU(5);\mathbb Z/2)$ is an isomorphism, the above equality gives
\begin{displaymath}
  w_{16}((\Delta_{10}^+)_{\mathbb R}) = y_6y_{10} + y_8^2 + y_4^4.
\end{displaymath}
Using the Wu formula, we obtain the total Stiefel--Whitney class
\begin{equation}
  \label{eq:delta10}
  \begin{split}
  w((\Delta_{10}^+)_{\mathbb R}) & = c(\Delta_{10}^+) \\
  & =
  1 + (y_6y_{10} + y_8^2 + y_4^4)
  \\ & \quad
  + (y_4y_{10}^2 + y_6y_8y_{10} + y_4^2y_6y_{10} + y_4^2y_8^2 + y_6^4)
  \\ & \quad
  + (y_8y_{10}^2 + y_4^2y_{10}^2 + y_4y_6y_8y_{10} + y_6^3y_{10}
  + y_6^2y_8^2 + y_7^4)
  \\ & \quad
  + (y_{10}^3 + y_4y_6y_{10}^2 + y_6^2y_8y_{10} + y_7^2y_8^2)
  + u_{32}.
  \end{split}
\end{equation}
Note that $c(\Delta_{10}^+) = c(\Delta_{10}^-)$, since $\Delta_{10}^+$
and $\Delta_{10}^-$ are conjugate to each other.
Using it, we can obtain the total Chern class $c(\rho_6)$ for the
representation $\rho_6 \co E_6 \to SU(27)$ mentioned in the
introduction, since we have by Corollary 8.3 of \cite{adams2} that the induced
representation
\begin{displaymath}
  \Spin(10) \longrightarrow E_6 \longrightarrow SU(27)
\end{displaymath}
is a direct sum of the one dimensional trivial representation,
$\Delta_{10}^+$ and the composition map
\begin{displaymath}
  \Spin(10) \longrightarrow SO(10) \longrightarrow SU(10).
\end{displaymath}
Recall from \eqref{eq:spin_rep} that $f_{10}^*((\Delta_{11})_{\mathbb C})
= \smash{\Delta_{10}^+} \oplus \smash{\Delta_{10}^-}$,
which gives $w(\smash{f_{10}^*}(\Delta_{11})_{\mathbb R}) =
w((\smash{\Delta_{10}^+})_{\mathbb R})w((\smash{\Delta_{10}^-})_{\mathbb R})$.
We consider the homomorphism
\begin{displaymath}
  Bf_{10}^* \co H^*(B\Spin(11);\mathbb Z/2) \longrightarrow
  H^*(B\Spin(10);\mathbb Z/2),
\end{displaymath}
where we see that a basis of $\mathrm{Ker} \; Bf_{10}^*$ of degree 32
is given by
\begin{displaymath}
  \{ y_{10}y_{11}^2,\
  y_4y_6y_{11}^2,\
  y_6y_7y_8y_{11},\
  y_7^3y_{11},\
  y_4^2y_6y_7y_{11}
  \}.
\end{displaymath}
So we can write as follows:
\begin{align*}
  w_{32}((\Delta_{11})_{\mathbb R}) & =
  a_1 y_{10}y_{11}^2
  + a_2 y_4y_6y_{11}^2
  + a_3 y_6y_7y_8y_{11}
  + a_4 y_7^3y_{11}
 \\ & \quad
  + a_5 y_4^2y_6y_7y_{11}
  + y_6^2y_{10}^2 + y_8^4 + y_4^8,
\end{align*}
where $a_i \in \mathbb Z/2$.
By \fullref{thm:quillen}, we have
\begin{displaymath}
  0 =
  \Sq^1 w_{32}((\Delta_{11})_{\mathbb R}) =
  a_1 y_{11}^3
  + a_2 y_4y_7y_{11}^2
  + a_3 y_7^2y_8y_{11}
  + a_5 y_4^2y_7^2y_{11},
\end{displaymath}
which implies $a_1 = a_2 = a_3$ and $a_5 = 0$.
Again by \fullref{thm:quillen}, we have
\begin{displaymath}
  0 =
  \Sq^4 w_{32}((\Delta_{11})_{\mathbb R}) =
  a_4 y_7^2y_{11}^2 + y_7^2y_{11}^2,
\end{displaymath}
which implies $a_4 = 1$.
Further by \fullref{thm:quillen}, we have
\begin{displaymath}
  0 = \Sq^{30} w_{32}((\Delta_{11})_{\mathbb R}) =
  (a_1 + 1)
  (y_4^2y_6^2y_{10}^2y_{11}^2
  + y_6^4y_8^2y_{11}^2
  + y_{10}^4y_{11}^2),
\end{displaymath}
which implies $a_1 = 1$.
Thus we obtain the total Stiefel--Whitney class
\begin{align*}
  w((\Delta_{11})_{\mathbb R}) & =
  c((\Delta_{11})_{\mathbb C}) = p(\Delta_{11}) \\
  & =
  1
  + (y_{10}y_{11}^2
  + y_4y_6y_{11}^2
  + y_6y_7y_8y_{11}
  + y_7^3y_{11}
  + y_6^2y_{10}^2 + y_8^4 + y_4^8)
  \\ & \quad
  + (y_8^2y_{10}y_{11}^2
  + y_4y_6y_8^2y_{11}^2
  + y_4^4y_{10}y_{11}^2
  + y_4^5y_6y_{11}^2
  + y_4y_6^2y_{10}y_{11}^2
  \\ & \qquad
  + y_6^3y_8y_{11}^2
  + y_6y_7y_8^3y_{11}
  + y_6y_{10}^2y_{11}^2
  + y_4^4y_6y_7y_8y_{11}
  \\ & \qquad
  + y_4^2y_7^3y_8y_{11}
  + y_4^3y_7^2y_{11}^2
  + y_4^4y_7^3y_{11}
  + y_7^3y_8^2y_{11}
  + y_6^2y_8^2y_{10}^2
  \\ & \qquad
  + y_4^2y_{10}^4
  + y_4^4y_6^2y_{10}^2
  + y_4^4y_8^4
  + y_6^8)
  \\ & \quad
  + (y_6^3y_8^2y_{11}^2
  + y_4^2y_8^2y_{10}y_{11}^2
  + y_4^3y_6y_8^2y_{11}^2
  + y_4y_6^5y_{11}^2
  + y_6^4y_{10}y_{11}^2
  \\ & \qquad
  + y_4y_{10}^3y_{11}^2
  + y_4^3y_6^2y_{10}y_{11}^2
  + y_4^2y_6^3y_8y_{11}^2
  + y_6y_8y_{10}^2y_{11}^2
  \\ & \qquad
  + y_4^2y_6y_7y_8^3y_{11}
  + y_6^5y_7y_8y_{11}
  + y_6^4y_7^3y_{11}
  + y_4^2y_6^2y_8^2y_{10}^2
  + y_6^6y_{10}^2
  \\ & \qquad
  + y_8^2y_{10}^4
  + y_4^4y_{10}^4
  + y_6^4y_8^4
  + y_7^8)
  \\ & \quad
  + (y_6^3y_{10}^2y_{11}^2
  + y_4^2y_{10}^3y_{11}^2
  + y_4^3y_6y_{10}^2y_{11}^2
  + y_4y_6^4y_{10}y_{11}^2
  + y_6^5y_8y_{11}^2
  \\ & \qquad
  + y_8y_{10}^3y_{11}^2
  + y_6^3y_7y_8^3y_{11}
  + y_6^4y_7^2y_{11}^2
  + y_7^7y_{11}
  + y_4^2y_7^2y_8^2y_{11}^2
  \\ & \qquad
  + y_4^2y_6^2y_{10}^4
  + y_6^4y_8^2y_{10}^2
  + y_8^2y_{11}^4
  + y_{10}^6
  + y_4^4y_{11}^4
  + y_7^4y_8^4)
  \\ & \quad
  + u_{64}.
\end{align*}
Recall from \eqref{eq:plus} that $\tilde{r}_6^*(\Delta_{12}^+)_{\mathbb C}
= 2 \oplus \lambda^2 \oplus \lambda^4$,
where $\lambda^i$ is the $i$--th exterior representation $SU(6) \to
SU(\smash{\tbinom{6}{i}})$.
In a similar way to the case of the $i$--th exterior representation
$SU(5) \to SU(\tbinom{5}{i})$, we obtain the mod 2 total Chern classes
\begin{align*}
  c(\lambda^2) = c(\lambda^4) & =
  1
  + (c_3c_5
  + c_4^2
  + c_2^4)
  \\ & \quad
  + (c_2c_5^2
  + c_3^2c_6
  + c_3c_4c_5
  + c_2^2c_3c_5
  + c_2^2c_4^2
  + c_3^4
  + c_6^2)
  \\ & \quad
  + (c_3c_5c_6
  + c_4c_5^2
  + c_2^2c_5^2
  + c_2c_3^2c_6
  + c_2c_3c_4c_5
  + c_3^3c_5
  + c_3^2c_4^2)
  \\ & \quad
  + (c_5^3
  + c_2c_3c_5^2
  + c_3^3c_6
  + c_3^2c_4c_5),
\end{align*}
where $c_i \in H^*(BSU(6);\mathbb Z/2)$ is the mod 2 $i$--th Chern class.
Then by \eqref{eq:whitney_sum},
\begin{displaymath}
  w_{32}(\tilde{r}_6^*(\Delta_{12}^+)_{\mathbb R})
  = c_{16}(\tilde{r}_6(\Delta_{12}^+)_{\mathbb C})
  = c_3^2c_5^2 + c_4^4 + c_2^8.
\end{displaymath}
Recall from \eqref{eq:spin_rep} that $f_{11}^*\Delta_{12}^+
= \Delta_{11}$,
which gives $f_{11}^*w((\Delta_{12}^+)_{\mathbb R}) =
w((\Delta_{11})_{\mathbb R})$.
We consider the homomorphism
\begin{displaymath}
  Bf_{11}^* \oplus B\tilde{r}_6^* \co H^*(B\Spin(12);\mathbb Z/2)
  \rightarrow H^*(B\Spin(11);\mathbb Z/2) \oplus
  H^*(BSU(6);\mathbb Z/2),
\end{displaymath}
where we see that a basis of $\mathrm{Ker} \; Bf_{11}^* \oplus B\tilde{r}_6^*$ of
degree 32 is given by
\begin{displaymath}
  \{ y_6y_7^2y_{12} \}.
\end{displaymath}
So we can write as follows:
\begin{align*}
  w_{32}((\Delta_{12}^+)_{\mathbb R}) =
  a_1 y_6y_7^2y_{12}
  &+ y_{10}y_{11}^2
  + y_4y_6y_{11}^2
  + y_6y_7y_8y_{11}
\\  &+ y_7^3y_{11}
  + y_6^2y_{10}^2
  + y_8^4
  + y_4^8,
\end{align*}
where $a_1 \in \mathbb Z/2$.
By \fullref{thm:quillen} we have
\begin{displaymath}
  0 = \Sq^1 w_{32}((\Delta_{12}^+)_{\mathbb R}) =
  a_1 y_7^3y_{12}
  + y_{11}^3
  + y_4y_7y_{11}^2
  + y_7^2y_8y_{11},
\end{displaymath}
which implies $a_1 = 1$.
Thus we have
\begin{displaymath}
  w_{32}((\Delta_{12}^+)_{\mathbb R}) {=}
  y_6y_7^2y_{12}
  + y_{10}y_{11}^2
  + y_4y_6y_{11}^2
  + y_6y_7y_8y_{11}
  + y_7^3y_{11}
  + y_6^2y_{10}^2
  + y_8^4
  + y_4^8,
\end{displaymath}
and the total Stiefel--Whitney class is given by
\begin{equation}
  \label{eq:delta12+}
  \begin{split}
  w((\Delta_{12}^+)_{\mathbb R}) & =
  c((\Delta_{12}^+)_{\mathbb C}) = p(\Delta_{12}^+) \\
  & =
  1 + w_{32}((\Delta_{12}^+)_{\mathbb R})
  + \Sq^{16}w_{32}((\Delta_{12}^+)_{\mathbb R})
  + \Sq^{24}w_{32}((\Delta_{12}^+)_{\mathbb R})
  \\ & \quad
  + \Sq^{28}w_{32}((\Delta_{12}^+)_{\mathbb R})
  + u_{64}.
  \end{split}
\end{equation}
In a similar way, we can obtain $w_{32}((\Delta_{12}^-)_{\mathbb R})$ which is
equal to $w_{32}((\Delta_{12}^+)_{\mathbb R})$.
In order to determine the total Stiefel--Whitney class of
$\smash{(\Delta_{12}^-)_{\mathbb R}}$, we need to calculate
$w_{64}(\smash{(\Delta_{12}^-)_{\mathbb R}})$.
From \eqref{eq:minus} we have 
$\smash{\tilde{r}_6^*((\Delta_{12}^-)_{\mathbb C})}
= \lambda^1 \oplus \lambda^3 \oplus \lambda^5$, for $\lambda^i$ the
$i$--th exterior representation
$SU(6) \to SU(\smash{\tbinom{6}{i}})$.
Similarly to the case of the $i$--th exterior representation
$SU(5) \to SU(\smash{\tbinom{5}{i}})$, we obtain the mod 2
total Chern class
\begin{align*}
  c(\lambda^3) & =
  1
  + c_2^2
  + c_3^2
  + (c_4^2
  + c_2^4)
  + c_5^2
  + (c_6^2
  + c_3^4
  + c_2^6)
  + c_2^4c_3^2
  \\ & \quad
  + (c_3^2c_5^2
  + c_2^4c_4^2
  + c_2^2c_3^4)
  + (c_2^4c_5^2
  + c_3^6)
  + (c_5^4
  + c_2^4c_6^2
  + c_2^2c_3^2c_5^2
  + c_3^4c_4^2), \\
  c(\lambda^1) & = c(\lambda^5)
  = 1 + c_2 + c_3 + c_4 + c_5 + c_6.
\end{align*}
Then by \eqref{eq:whitney_sum} we obtain
\begin{displaymath}
  w_{64}(\tilde{r}_6^*(\Delta_{12}^-)_{\mathbb R}) = c_{32}(\tilde{r}_6^*(\Delta_{12}^-
)_{\mathbb C}) =
  c_6^2
  (c_5^4
  + c_2^4c_6^2
  + c_2^2c_3^2c_5^2
  + c_3^4c_4^2).
\end{displaymath}
Recall from \eqref{eq:spin_rep} that $f_{11}^*\Delta_{12}^-
= \Delta_{11}$,
which gives $w(f_{11}^*(\Delta_{12}^-)_{\mathbb R}) =
w((\Delta_{11})_{\mathbb R})$,
and hence we obtain
\begin{align*}
  w_{64}((\Delta_{12}^-)_{\mathbb R}) &\equiv u_{64}
  + y_{10}^4y_{12}^2
  + y_4^4y_{12}^4
  + y_4^2y_6^2y_{10}^2y_{12}^2
  + y_6^4y_8^2y_{12}^2\\
  &\qquad\mod \mathrm{Ker} \; Bf_{11}^* \oplus B\tilde{r}_6^*,\\
\tag*{\hbox{where}}
  Bf_{11}^* \co &H^*(B\Spin(12);\mathbb Z/2) \longrightarrow
  H^*(B\Spin(11);\mathbb Z/2), \\
  B\tilde{r}_6^* \co &H^*(B\Spin(12);\mathbb Z/2) \longrightarrow
  H^*(BSU(6);\mathbb Z/2).
\end{align*}
In order to determine $w_{64}((\Delta_{12}^-)_{\mathbb R})$ exactly, we use the
method of indeterminate coefficients;
using the equations $\Sq^nw_{64}((\Delta_{12}^-)_{\mathbb R} = 0$ for
$n = 1, 2, 4, 62$,
we can determine all the indeterminate coefficients,
and hence the 64--th Stiefel--Whitney class is given as follows:
\begin{align*}
  w_{64}((\Delta_{12}^-)_{\mathbb R}) & = u_{64}
  + y_{10}^4y_{12}^2
  + y_{10}^3y_{11}^2y_{12}
  + y_7^4y_{12}^3
  + y_4^2y_{10}y_{11}^2y_{12}^2
    \\ & \quad
  + y_4y_7^2y_{11}^2y_{12}^2
  + y_6^3y_{11}^2y_{12}^2
  + y_7^3y_8y_{11}y_{12}^2
  + y_7^2y_8^2y_{11}^2y_{12}
  \\ & \quad
  + y_4^4y_{12}^4
  + y_4^2y_6y_7^2y_{12}^3
  + y_4^3y_6y_{11}^2y_{12}^2
  + y_4^2y_6y_7y_8y_{11}y_{12}^2
 \\ & \quad
  + y_4^2y_7^3y_{11}y_{12}^2
  + y_4y_6^3y_7y_{11}y_{12}^2
  + y_4^2y_6^2y_{10}^2y_{12}^2
  + y_6^4y_8^2y_{12}^2
    \\ & \quad
  + y_6^3y_7^2y_8y_{12}^2
  + y_4^2y_6^2y_{10}y_{11}^2y_{12}
  + y_4y_7^3y_8^2y_{11}y_{12}
    \\ & \quad
  + y_6^3y_7y_8^2y_{11}y_{12}
  + y_7^4y_8^3y_{12}.
\end{align*}
Thus we obtain the total Stiefel--Whitney class
\begin{equation}
  \label{eq:delta12-}
  \begin{split}
  w((\Delta_{12}^-)_{\mathbb R}) & =
  c((\Delta_{12}^-)_{\mathbb C}) = p(\Delta_{12}^-) \\
  & =
  1 + w_{32}((\Delta_{12}^+)_{\mathbb R})
  + \Sq^{16}w_{32}((\Delta_{12}^+)_{\mathbb R})
  + \Sq^{24}w_{32}((\Delta_{12}^+)_{\mathbb R})
  \\ & \quad
  + \Sq^{28}w_{32}((\Delta_{12}^+)_{\mathbb R})
  + w_{64}((\Delta_{12}^-)_{\mathbb R}).
  \end{split}
\end{equation}
Using it, we can obtain the total Pontrjagin class $p(\rho_7)$ for $\rho_7 \co E_7 \to Sp(28)$,
the representation  mentioned in the
introduction, since we have by Corollary 8.2 of \cite{adams2} that the induced
representation
\begin{displaymath}
  \Spin(12) \longrightarrow E_7 \longrightarrow Sp(28)
\end{displaymath}
is a direct sum of
\begin{displaymath}
  \Delta_{12}^- \co \Spin(12) \longrightarrow Sp(16)
\end{displaymath}
and the composition map
\begin{displaymath}
  \Spin(12) \longrightarrow SO(12) \longrightarrow SU(12)
  \longrightarrow Sp(12).
\end{displaymath}
Recall from \eqref{eq:spin_rep} that
$f_{12}^*\Delta_{13} = \Delta_{12}^+ \oplus \Delta_{12}^-$.
Then we have $$w(f_{12}^*(\Delta_{13})_{\mathbb R} =
w((\Delta_{12}^+)_{\mathbb R})w((\Delta_{12}^-)_{\mathbb R}).$$
In order to determine $w_{64}((\Delta_{13})_{\mathbb R})$, we use the
method of indeterminate coefficients;
using the equations $\Sq^nw_{64}((\Delta_{13})_{\mathbb R}) = 0$ for
$n = 1, 2, 4, 8, 62$,
we can determine all the indeterminate coefficients,
and hence the 64--th Stiefel--Whitney class is
given as follows:
\begin{align*}
  w_{64}((\Delta&_{13})_{\mathbb R})  =
  y_{12}y_{13}^4
  + y_6^2y_{13}^4
  + y_4y_{10}y_{11}y_{13}^3
  + y_7^2y_{11}y_{13}^3
  + y_4y_{11}^2y_{12}y_{13}^2
  \\ & \quad
  + y_7y_8y_{11}y_{12}y_{13}^2
  + y_6y_{10}^2y_{12}y_{13}^2
  + y_6y_{10}y_{11}^2y_{13}^2
  + y_8y_{10}^2y_{11}y_{12}y_{13}
  \\ & \quad
  + y_{10}^4y_{11}y_{13}
  + y_4^3y_{13}^4
  + y_4^2y_6y_{11}y_{13}^3
  + y_4y_7^3y_{13}^3
  + y_6^3y_8y_{12}y_{13}^2
  \\ & \quad
  + y_6^2y_7^2y_{12}y_{13}^2
  + y_6^2y_7y_8y_{11}y_{13}^2
  + y_6y_7^3y_{11}y_{13}^2
  + y_4^2y_8y_{11}y_{12}^2y_{13}
  \\ & \quad
  + y_4y_6y_7y_{11}^2y_{12}y_{13}
  + y_4^2y_{10}^2y_{11}y_{12}y_{13}
  + y_6^3y_{10}y_{11}y_{12}y_{13}
  \\ & \quad
  + y_6^2y_8^2y_{11}y_{12}y_{13}
  + y_6y_7^2y_8y_{11}y_{12}y_{13}
  + y_7^4y_{11}y_{12}y_{13}
  + y_7^3y_8y_{11}^2y_{13}
  \\ & \quad
  + y_4^3y_6y_7y_{13}^3
  + y_4^3y_7^2y_{12}y_{13}^2
  + y_4^2y_6^3y_{12}y_{13}^2
  + y_4^2y_6^2y_7y_{11}y_{13}^2
  \\ & \quad
  + y_4^2y_7^2y_8^2y_{13}^2
  + y_4y_6^2y_7^2y_8y_{13}^2
  + y_4y_6y_7^4y_{13}^2
  + y_4^3y_6y_{10}y_{11}y_{12}y_{13}
  \\ & \quad
  + y_4^2y_6y_7^2y_{11}y_{12}y_{13}
  + y_4y_6^4y_{11}y_{12}y_{13}
  + y_6^4y_7y_8y_{12}y_{13}
  + y_4y_7^5y_{12}y_{13}
  \\ & \quad
  + y_4^2y_6^2y_{10}^2y_{11}y_{13}
  + y_6^4y_8^2y_{11}y_{13}
  + y_4y_7^4y_8y_{11}y_{13}
  + y_7^5y_8^2y_{13}
  + y_{10}^4y_{12}^2
  \\ & \quad
  + y_{10}^3y_{11}^2y_{12}
  + y_7^4y_{12}^3
  + y_4^2y_{10}y_{11}^2y_{12}^2
  + y_4y_7^2y_{11}^2y_{12}^2
  + y_6^3y_{11}^2y_{12}^2
  \\ & \quad
  + y_7^3y_8y_{11}y_{12}^2
  + y_7^2y_8^2y_{11}^2y_{12}
  + y_4^4y_{12}^4
  + y_4^2y_6y_7^2y_{12}^3
  + y_4^3y_6y_{11}^2y_{12}^2
  \\ & \quad
  + y_4^2y_6y_7y_8y_{11}y_{12}^2
  + y_4^2y_7^3y_{11}y_{12}^2
  + y_4y_6^3y_7y_{11}y_{12}^2
  + y_4^2y_6^2y_{10}^2y_{12}^2
  \\ & \quad
  + y_6^4y_8^2y_{12}^2
  + y_6^3y_7^2y_8y_{12}^2
  + y_4^2y_6^2y_{10}y_{11}^2y_{12}
  + y_4y_7^3y_8^2y_{11}y_{12}
  \\ & \quad
  + y_6^3y_7y_8^2y_{11}y_{12}
  + y_7^4y_8^3y_{12}
  + y_6^2y_7^4y_{12}^2
  + y_{10}^2y_{11}^4
  + y_4^2y_6^2y_{11}^4
  \\ & \quad
  + y_6^2y_7^2y_8^2y_{11}^2
  + y_7^6y_{11}^2
  + y_6^4y_{10}^4
  + y_8^8
  + y_4^{16}
  .
\end{align*}
Thus we obtain the total Stiefel--Whitney class
\begin{align*}
  w((\Delta_{13})_{\mathbb R}) & =
  c((\Delta_{13})_{\mathbb C}) = p(\Delta_{13}) \\
  & =
  1 + w_{64}((\Delta_{13})_{\mathbb R})
  + \Sq^{32}w_{64}((\Delta_{13})_{\mathbb R})
  + \Sq^{48}w_{64}((\Delta_{13})_{\mathbb R})
  \\ & \quad
  + \Sq^{56}w_{64}((\Delta_{13})_{\mathbb R})
  + \Sq^{60}w_{64}((\Delta_{13})_{\mathbb R})
  + u_{128}.
\end{align*}
Recall from \eqref{eq:plus} that if $\lambda^i$ is the $i$--th exterior representation $SU(7) \to
SU(\smash{\tbinom{7}{i}})$ then
$\smash{\tilde{r}_7^*((\Delta_{14}^+)_{\mathbb C})} = 1 \oplus \lambda^2
\oplus \lambda^4 \oplus \lambda^6$.
In a similar way to the case of the $i$--th exterior representation
$SU(5) \to SU(\smash{\tbinom{5}{i}})$,
we obtain the mod 2 total Chern classes
\begin{align*}
  c(\lambda^2) & =
  1
 {+}c_2
 {+}c_3
 {+}c_4
 {+}c_5
 {+}c_6
 {+}c_7
 {+}(c_3c_5
 {+}c_4^2
 {+}c_2^4)
 {+}(c_2c_3c_5
 {+}c_2c_4^2
 {+}c_2^5)
  \\ & \quad
 {+}(c_3^2c_5
 {+}c_3c_4^2
 {+}c_2^4c_3)
  \\ & \quad
 {+}(c_6^2
 {+}c_2c_5^2
 {+}c_3^2c_6
 {+}c_4^3
 {+}c_2^2c_3c_5
 {+}c_2^2c_4^2
 {+}c_3^4
 {+}c_2^4c_4
 {+}c_2c_3c_7)
  \\ & \quad
 {+}(c_3c_5^2
 {+}c_4^2c_5
 {+}c_2^4c_5)
  \\ & \quad
 {+}(c_2c_6^2
 {+}c_4^2c_6
 {+}c_4c_5^2
 {+}c_3^3c_5
 {+}c_3^2c_4^2
 {+}c_2^4c_6
 {+}c_2^3c_3c_5
 {+}c_2^3c_4^2
 {+}c_2c_3^4
  \\ & \qquad
 {+}c_2c_5c_7
 {+}c_7^2)
  \\ & \quad
 {+}(c_3c_6^2
 {+}c_5^3
 {+}c_2^2c_3^2c_5
 {+}c_2^2c_3c_4^2
 {+}c_3^5
 {+}c_3c_5c_7
 {+}c_4^2c_7
 {+}c_2^4c_7)
  \\ & \quad
 {+}(c_4c_6^2
 {+}c_2c_3c_5c_6
 {+}c_3^2c_4c_6
 {+}c_3c_4^2c_5
 {+}c_2^3c_5^2
 {+}c_2^2c_3^2c_6
 {+}c_2^2c_4^3
 {+}c_2c_3^3c_5
  \\ & \qquad
 {+}c_2c_3^2c_4^2
 {+}c_3^4c_4
 {+}c_2c_3c_4c_7
 {+}c_3^3c_7
 {+}c_2^3c_3c_7
 {+}c_2c_7^2)
  \\ & \quad
 {+}(c_5c_6^2
 {+}c_2^2c_3c_5^2
 {+}c_2^2c_4^2c_5
 {+}c_3^3c_4^2
 {+}c_3c_7^2)
  \\ & \quad
 {+}(c_6^3
 {+}c_2c_5^2c_6
 {+}c_3^2c_6^2
 {+}c_3c_5^3
 {+}c_4^2c_5^2
 {+}c_2^2c_3c_5c_6
 {+}c_2^2c_4^2c_6
 {+}c_2^2c_4c_5^2
  \\ & \qquad
 {+}c_2c_3^2c_4c_6
 {+}c_2c_3^2c_5^2
 {+}c_2c_3c_4^2c_5
 {+}c_3^2c_4^3
 {+}c_2c_3c_6c_7
 {+}c_2c_4c_5c_7
 {+}c_2^3c_5c_7
  \\ & \qquad
 {+}c_2^2c_3c_4c_7
 {+}c_4c_7^2)
  \\ & \quad
 {+}(c_3c_5^2c_6
 {+}c_2^2c_5^3
 {+}c_2c_3^2c_5c_6
 {+}c_3^3c_4c_6
 {+}c_3^3c_5^2
 {+}c_6^2c_7
 {+}c_3^2c_6c_7
 {+}c_3c_4c_5c_7
  \\ & \qquad
 {+}c_2^2c_3c_5c_7
 {+}c_2^2c_4^2c_7
 {+}c_2c_3^2c_4c_7
 {+}c_5c_7^2
 {+}c_2c_3c_7^2)
  \\ & \quad
 {+}(c_3c_5c_6^2
 {+}c_4c_5^2c_6
 {+}c_5^4
 {+}c_2^2c_5^2c_6
 {+}c_2c_3^2c_6^2
 {+}c_2c_3c_4c_5c_6
 {+}c_2c_3c_5^3
 {+}c_3^2c_4^2c_6
  \\ & \qquad
 {+}c_3^2c_4c_5^2
 {+}c_2c_5c_6c_7
 {+}c_3c_5^2c_7
 {+}c_2^2c_4c_5c_7
 {+}c_6c_7^2
 {+}c_2^3c_7^2)
  \\ & \quad
 {+}(c_5^3c_6
 {+}c_2c_3c_5^2c_6
 {+}c_3^3c_6^2
 {+}c_3^2c_4c_5c_6
 {+}c_3c_5c_6c_7
 {+}c_4c_5^2c_7
 {+}c_2c_3c_4c_5c_7
  \\ & \qquad
 {+}c_3^2c_4^2c_7
 {+}c_2c_5c_7^2
 {+}c_2^2c_3c_7^2
 {+}c_7^3)
  , \\
  c(\lambda^4) & =
  1
 {+}c_2^2
 {+}c_3^2
 {+}c_3c_5
 {+}c_5^2
 {+}(c_2c_5^2
 {+}c_3^2c_6
 {+}c_3c_4c_5
 {+}c_2c_3c_7)
  \\ & \quad
 {+}(c_3c_5c_6
 {+}c_4c_5^2
 {+}c_2^2c_5^2
 {+}c_2c_3^2c_6
 {+}c_2c_3c_4c_5
 {+}c_2c_5c_7
 {+}c_2^2c_3c_7)
  \\ & \quad
 {+}(c_5^3
 {+}c_2c_3c_5^2
 {+}c_3^3c_6
 {+}c_3^2c_4c_5
 {+}c_2c_3^2c_7)
  \\ & \quad
 {+}(c_2^2c_6^2
 {+}c_3^2c_5^2
 {+}c_3c_4^2c_5
 {+}c_4^4
 {+}c_2^3c_5^2
 {+}c_2^2c_3^2c_6
 {+}c_2^2c_3c_4c_5
 {+}c_2^8
 {+}c_2^2c_5c_7
  \\ & \qquad
 {+}c_3^3c_7
 {+}c_2^3c_3c_7)
  \\ & \quad
 {+}(c_3^2c_6^2
 {+}c_3c_5^3
 {+}c_4^2c_5^2
 {+}c_2^2c_3c_5c_6
 {+}c_2^2c_4c_5^2
 {+}c_2c_3^2c_5^2
 {+}c_3^4c_6
 {+}c_3^3c_4c_5
  \\ & \qquad
 {+}c_2^4c_5^2
 {+}c_2^3c_3^2c_6
 {+}c_2^3c_3c_4c_5
 {+}c_2^3c_5c_7
 {+}c_2c_3^3c_7
 {+}c_2^4c_3c_7
 {+}c_2^2c_7^2)
  \\ & \quad
 {+}(c_2^2c_5^3
 {+}c_2^3c_3c_5^2
 {+}c_2^2c_3^3c_6
 {+}c_2^2c_3^2c_4c_5
 {+}c_2^3c_3^2c_7)
  \\ & \quad
 {+}(c_3c_5c_6^2
 {+}c_5^4
 {+}c_2c_4^2c_5^2
 {+}c_3^3c_5c_6
 {+}c_3^2c_4^2c_6
 {+}c_3^2c_4c_5^2
 {+}c_3c_4^3c_5
 {+}c_2^4c_6^2
  \\ & \qquad
 {+}c_2^2c_3c_4^2c_5
 {+}c_2^2c_4^4
 {+}c_2c_3^4c_6
 {+}c_2c_3^3c_4c_5
 {+}c_2^5c_5^2
 {+}c_2^4c_3^2c_6
 {+}c_2^4c_3c_4c_5
 {+}c_2^{10}
  \\ & \qquad
 {+}c_3c_5^2c_7
 {+}c_2^2c_3c_6c_7
 {+}c_2c_3c_4^2c_7
 {+}c_3^3c_4c_7
 {+}c_2^4c_5c_7
 {+}c_2^5c_3c_7
 {+}c_3^2c_7^2
  \\ & \qquad
 {+}c_2^3c_7^2)
  \\ & \quad
 {+}(c_3^2c_5^3
 {+}c_2c_3^3c_5^2
 {+}c_3^5c_6
 {+}c_3^4c_4c_5
 {+}c_2c_3^4c_7)
  \\ & \quad
 {+}(c_5^2c_6^2
 {+}c_2c_5^4
 {+}c_3^2c_5^2c_6
 {+}c_3c_4^2c_5c_6
 {+}c_3c_4c_5^3
 {+}c_4^3c_5^2
 {+}c_2^2c_3c_5^3
 {+}c_2c_3^2c_4^2c_6
  \\ & \qquad
 {+}c_2c_3c_4^3c_5
 {+}c_3^4c_5^2
 {+}c_3^3c_4^2c_5
 {+}c_3^2c_4^4
 {+}c_2^4c_3c_5c_6
 {+}c_2^4c_4c_5^2
 {+}c_2^6c_5^2
 {+}c_2^5c_3^2c_6
  \\ & \qquad
 {+}c_2^5c_3c_4c_5
 {+}c_2^8c_3^2
 {+}c_5^3c_7
 {+}c_2^2c_5c_6c_7
 {+}c_2c_3c_5^2c_7
 {+}c_2c_4^2c_5c_7
 {+}c_3^2c_4c_5c_7
  \\ & \qquad
 {+}c_2^3c_3c_6c_7
 {+}c_2^2c_3c_4^2c_7
 {+}c_2c_3^3c_4c_7
 {+}c_3^5c_7
 {+}c_2^5c_5c_7
 {+}c_2^6c_3c_7
 {+}c_3c_5c_7^2
  \\ & \qquad
 {+}c_2c_3^2c_7^2)
  \\ & \quad
 {+}(c_4^2c_5^3
 {+}c_2c_3c_4^2c_5^2
 {+}c_3^3c_4^2c_6
 {+}c_3^2c_4^3c_5
 {+}c_2^4c_5^3
 {+}c_2^5c_3c_5^2
 {+}c_2^4c_3^3c_6
 {+}c_2^4c_3^2c_4c_5
  \\ & \qquad
 {+}c_3^2c_5^2c_7
 {+}c_2^2c_3^2c_6c_7
 {+}c_2c_3^3c_5c_7
 {+}c_2c_3^2c_4^2c_7
 {+}c_3^4c_4c_7
 {+}c_2^5c_3^2c_7
 {+}c_2^2c_5c_7^2
  \\ & \qquad
 {+}c_3^3c_7^2
 {+}c_2^3c_3c_7^2)
  \\ & \quad
 {+}(c_6^4
 {+}c_2c_5^2c_6^2
 {+}c_3^2c_6^3
 {+}c_3c_4c_5c_6^2
 {+}c_3c_5^3c_6
 {+}c_4c_5^4
 {+}c_2^2c_3c_5c_6^2
 {+}c_2^2c_4^2c_6^2
  \\ & \qquad
 {+}c_2^2c_5^4
 {+}c_2c_3^2c_5^2c_6
 {+}c_2c_3c_4c_5^3
 {+}c_2c_3^4c_5^2
 {+}c_3^6c_6
 {+}c_3^5c_4c_5
 {+}c_2^6c_6^2
 {+}c_2^4c_3^2c_5^2
  \\ & \qquad
 {+}c_2^4c_3c_4^2c_5
 {+}c_2^4c_4^4
 {+}c_3^8
 {+}c_2^7c_5^2
 {+}c_2^6c_3^2c_6
 {+}c_2^6c_3c_4c_5
 {+}c_2^8c_3c_5
 {+}c_2c_3c_6^2c_7
  \\ & \qquad
 {+}c_3c_4c_5^2c_7
 {+}c_2^3c_5c_6c_7
 {+}c_2^2c_3c_4c_6c_7
 {+}c_2^2c_4^2c_5c_7
 {+}c_2c_3^5c_7
 {+}c_2^6c_5c_7
  \\ & \qquad
 {+}c_2^4c_3^3c_7
 {+}c_2^7c_3c_7
 {+}c_5^2c_7^2
 {+}c_2^3c_4c_7^2)
  \\ & \quad
 {+}(c_5^5
 {+}c_2c_3c_5^4
 {+}c_3^3c_5^2c_6
 {+}c_3^2c_4c_5^3
 {+}c_2c_3^2c_5^2c_7)
  \\ & \quad
 {+}(c_3c_5c_6^3
 {+}c_4c_5^2c_6^2
 {+}c_2^2c_5^2c_6^2
 {+}c_2c_3^2c_6^3
 {+}c_2c_3c_4c_5c_6^2
 {+}c_3^3c_5c_6^2
 {+}c_3^2c_4^2c_6^2
  \\ & \qquad
 {+}c_3^5c_5c_6
 {+}c_3^4c_4c_5^2
 {+}c_2^4c_3^2c_6^2
 {+}c_2^4c_3c_5^3
 {+}c_2^4c_4^2c_5^2
 {+}c_2^2c_3^4c_5^2
 {+}c_2c_3^6c_6
  \\ & \qquad
 {+}c_2c_3^5c_4c_5
 {+}c_2^6c_3c_5c_6
 {+}c_2^6c_4c_5^2
 {+}c_2^5c_3^2c_5^2
 {+}c_2^4c_3^4c_6
 {+}c_2^4c_3^3c_4c_5
 {+}c_2^7c_3^2c_6
  \\ & \qquad
 {+}c_2^7c_3c_4c_5
 {+}c_2c_5c_6^2c_7
 {+}c_3c_5^2c_6c_7
 {+}c_4c_5^3c_7
 {+}c_2^2c_4c_5c_6c_7
 {+}c_2^2c_5^3c_7
  \\ & \qquad
 {+}c_2c_3^2c_5c_6c_7
 {+}c_2c_3c_4c_5^2c_7
 {+}c_3^3c_4c_6c_7
 {+}c_3^3c_5^2c_7
 {+}c_3^2c_4^2c_5c_7
 {+}c_2c_3^4c_5c_7
  \\ & \qquad
 {+}c_2^2c_3^5c_7
 {+}c_2^7c_5c_7
 {+}c_2^5c_3^3c_7
 {+}c_2^8c_3c_7
 {+}c_2c_5^2c_7^2
 {+}c_3^2c_6c_7^2
 {+}c_3c_4c_5c_7^2
  \\ & \qquad
 {+}c_2^3c_6c_7^2
 {+}c_2^2c_4^2c_7^2
 {+}c_2c_3^2c_4c_7^2
 {+}c_2^6c_7^2
 {+}c_2c_3c_7^3)
  \\ & \quad
 {+}(c_5^3c_6^2
 {+}c_2c_3c_5^2c_6^2
 {+}c_3^3c_6^3
 {+}c_3^2c_4c_5c_6^2
 {+}c_3^4c_5^3
 {+}c_2c_3^5c_5^2
 {+}c_3^7c_6
 {+}c_3^6c_4c_5
  \\ & \qquad
 {+}c_2^6c_5^3
 {+}c_2^7c_3c_5^2
 {+}c_2^6c_3^3c_6
 {+}c_2^6c_3^2c_4c_5
 {+}c_5^4c_7
 {+}c_2^2c_5^2c_6c_7
 {+}c_2c_3^2c_6^2c_7
  \\ & \qquad
 {+}c_2c_3c_5^3c_7
 {+}c_3^2c_4c_5^2c_7
 {+}c_2c_3^6c_7
 {+}c_2^7c_3^2c_7
 {+}c_3c_5^2c_7^2
 {+}c_2^2c_3c_6c_7^2
  \\ & \qquad
 {+}c_2^2c_4c_5c_7^2
 {+}c_2^3c_7^3)
  \\ & \quad
 {+}(c_3^2c_5^2c_6^2
 {+}c_3c_5^5
 {+}c_4^2c_5^4
 {+}c_2c_3^2c_5^4
 {+}c_3^4c_5^2c_6
 {+}c_3^3c_4c_5^3
 {+}c_2^4c_3c_5c_6^2
 {+}c_2^4c_5^4
  \\ & \qquad
 {+}c_2^2c_3^4c_6^2
 {+}c_3^6c_5^2
 {+}c_3^5c_4^2c_5
 {+}c_3^4c_4^4
 {+}c_2^5c_4^2c_5^2
 {+}c_2^4c_3^3c_5c_6
 {+}c_2^4c_3^2c_4^2c_6
  \\ & \qquad
 {+}c_2^4c_3^2c_4c_5^2
 {+}c_2^4c_3c_4^3c_5
 {+}c_2^3c_3^4c_5^2
 {+}c_2^2c_3^6c_6
 {+}c_2^2c_3^5c_4c_5
 {+}c_2^8c_6^2
 {+}c_2^6c_3c_4^2c_5
  \\ & \qquad
 {+}c_2^6c_4^4
 {+}c_2^5c_3^4c_6
 {+}c_2^5c_3^3c_4c_5
 {+}c_2^2c_3^8
 {+}c_5^3c_6c_7
 {+}c_2c_3c_5^2c_6c_7
 {+}c_3^3c_6^2c_7
  \\ & \qquad
 {+}c_3^2c_4c_5c_6c_7
 {+}c_2c_3^3c_5^2c_7
 {+}c_2^4c_3c_5^2c_7
 {+}c_2^2c_3^4c_5c_7
 {+}c_3^7c_7
 {+}c_2^6c_3c_6c_7
  \\ & \qquad
 {+}c_2^5c_3c_4^2c_7
 {+}c_2^4c_3^3c_4c_7
 {+}c_2^3c_3^5c_7
 {+}c_2^8c_5c_7
 {+}c_3c_5c_6c_7^2
 {+}c_4c_5^2c_7^2
 {+}c_2^2c_5^2c_7^2
  \\ & \qquad
 {+}c_2c_3c_4c_5c_7^2
 {+}c_3^2c_4^2c_7^2
 {+}c_2^4c_3^2c_7^2
 {+}c_2^7c_7^2
 {+}c_2c_5c_7^3
 {+}c_2^2c_3c_7^3
 {+}c_7^4)
  \\ & \quad
 {+}(c_2^4c_3^2c_5^3
 {+}c_2^5c_3^3c_5^2
 {+}c_2^4c_3^5c_6
 {+}c_2^4c_3^4c_4c_5
 {+}c_2^5c_3^4c_7)
  \\ & \quad
 {+}(c_3^3c_5^3c_6
 {+}c_3^2c_4c_5^4
 {+}c_2^4c_5^2c_6^2
 {+}c_2^2c_3^2c_5^4
 {+}c_2c_3^4c_5^2c_6
 {+}c_2c_3^3c_4c_5^3
 {+}c_3^6c_6^2
 {+}c_3^5c_5^3
  \\ & \qquad
 {+}c_3^4c_4^2c_5^2
 {+}c_2^5c_5^4
 {+}c_2^4c_3^2c_5^2c_6
 {+}c_2^4c_3c_4^2c_5c_6
 {+}c_2^4c_3c_4c_5^3
 {+}c_2^4c_4^3c_5^2
  \\ & \qquad
 {+}c_2^2c_3^5c_5c_6
 {+}c_2^2c_3^4c_4c_5^2
 {+}c_2c_3^6c_5^2
 {+}c_3^8c_6
 {+}c_3^7c_4c_5
 {+}c_2^6c_3c_5^3
 {+}c_2^5c_3^2c_4^2c_6
  \\ & \qquad
 {+}c_2^5c_3c_4^3c_5
 {+}c_2^4c_3^3c_4^2c_5
 {+}c_2^4c_3^2c_4^4
 {+}c_2^3c_3^6c_6
 {+}c_2^3c_3^5c_4c_5
 {+}c_3^{10}
 {+}c_2c_3^2c_5^3c_7
  \\ & \qquad
 {+}c_2^4c_5^3c_7
 {+}c_2^2c_3^3c_5^2c_7
 {+}c_2^6c_5c_6c_7
 {+}c_2^5c_3c_5^2c_7
 {+}c_2^5c_4^2c_5c_7
 {+}c_2^4c_3^2c_4c_5c_7
  \\ & \qquad
 {+}c_2^3c_3^4c_5c_7
 {+}c_2c_3^7c_7
 {+}c_2^7c_3c_6c_7
 {+}c_2^6c_3c_4^2c_7
 {+}c_2^5c_3^3c_4c_7
 {+}c_2^4c_3c_5c_7^2
  \\ & \qquad
 {+}c_2^2c_3^4c_7^2
 {+}c_2^5c_3^2c_7^2)
  \\ & \quad
 {+}(c_3^2c_5^5
 {+}c_2c_3^3c_5^4
 {+}c_3^5c_5^2c_6
 {+}c_3^4c_4c_5^3
 {+}c_2^4c_4^2c_5^3
 {+}c_2^2c_3^4c_5^3
 {+}c_2^5c_3c_4^2c_5^2
  \\ & \qquad
 {+}c_2^4c_3^3c_4^2c_6
 {+}c_2^4c_3^2c_4^3c_5
 {+}c_2^3c_3^5c_5^2
 {+}c_2^2c_3^7c_6
 {+}c_2^2c_3^6c_4c_5
 {+}c_2c_3^4c_5^2c_7
  \\ & \qquad
 {+}c_2^4c_3^2c_5^2c_7
 {+}c_2^6c_3^2c_6c_7
 {+}c_2^5c_3^3c_5c_7
 {+}c_2^5c_3^2c_4^2c_7
 {+}c_2^4c_3^4c_4c_7
 {+}c_2^3c_3^6c_7
  \\ & \qquad
 {+}c_2^6c_5c_7^2
 {+}c_2^4c_3^3c_7^2
 {+}c_2^7c_3c_7^2)
  \\ & \quad
 {+}(c_5^4c_6^2
 {+}c_2c_5^6
 {+}c_3^2c_5^4c_6
 {+}c_3c_4c_5^5
 {+}c_2^4c_6^4
 {+}c_2^2c_3^2c_5^2c_6^2
 {+}c_2^2c_3c_5^5
 {+}c_2^2c_4^2c_5^4
  \\ & \qquad
 {+}c_3^4c_4^2c_6^2
 {+}c_3^4c_5^4
 {+}c_2^5c_5^2c_6^2
 {+}c_2^4c_3^2c_6^3
 {+}c_2^4c_3c_4c_5c_6^2
 {+}c_2^4c_3c_5^3c_6
 {+}c_2^4c_4c_5^4
  \\ & \qquad
 {+}c_2^3c_3^2c_5^4
 {+}c_2^2c_3^4c_5^2c_6
 {+}c_2^2c_3^3c_4c_5^3
 {+}c_2c_3^4c_4^2c_5^2
 {+}c_3^7c_5c_6
 {+}c_3^6c_4^2c_6
 {+}c_3^6c_4c_5^2
  \\ & \qquad
 {+}c_3^5c_4^3c_5
 {+}c_2^6c_3c_5c_6^2
 {+}c_2^6c_4^2c_6^2
 {+}c_2^6c_5^4
 {+}c_2^5c_3^2c_5^2c_6
 {+}c_2^5c_3c_4c_5^3
 {+}c_2^4c_3^4c_6^2
  \\ & \qquad
 {+}c_2^2c_3^5c_4^2c_5
 {+}c_2^2c_3^4c_4^4
 {+}c_2c_3^8c_6
 {+}c_2c_3^7c_4c_5
 {+}c_3^9c_5
 {+}c_5^5c_7
 {+}c_2^2c_5^3c_6c_7
  \\ & \qquad
 {+}c_2c_3c_5^4c_7
 {+}c_3^3c_5^2c_6c_7
 {+}c_2^3c_3c_5^2c_6c_7
 {+}c_2^2c_3^3c_6^2c_7
 {+}c_2^2c_3^2c_4c_5c_6c_7
  \\ & \qquad
 {+}c_2c_3^4c_5c_6c_7
 {+}c_3^5c_4c_6c_7
 {+}c_3^4c_4^2c_5c_7
 {+}c_2^5c_3c_6^2c_7
 {+}c_2^4c_3c_4c_5^2c_7
 {+}c_2^3c_3^3c_5^2c_7
  \\ & \qquad
 {+}c_2^2c_3^5c_6c_7
 {+}c_2c_3^5c_4^2c_7
 {+}c_3^7c_4c_7
 {+}c_2^7c_5c_6c_7
 {+}c_2^6c_3c_4c_6c_7
 {+}c_2^6c_4^2c_5c_7
  \\ & \qquad
 {+}c_2^4c_3^4c_5c_7
 {+}c_3c_5^3c_7^2
 {+}c_2^2c_3c_5c_6c_7^2
 {+}c_2^2c_4c_5^2c_7^2
 {+}c_3^4c_6c_7^2
 {+}c_3^3c_4c_5c_7^2
  \\ & \qquad
 {+}c_2^3c_3c_4c_5c_7^2
 {+}c_2^2c_3^2c_4^2c_7^2
 {+}c_2c_3^4c_4c_7^2
 {+}c_3^6c_7^2
 {+}c_2^3c_3^4c_7^2
 {+}c_2^7c_4c_7^2
 {+}c_2^3c_5c_7^3
  \\ & \qquad
 {+}c_2c_3^3c_7^3
 {+}c_2^4c_3c_7^3)
  \\ & \quad
 {+}(c_2^4c_5^5
 {+}c_3^6c_5^3
 {+}c_2^5c_3c_5^4
 {+}c_2^4c_3^3c_5^2c_6
 {+}c_2^4c_3^2c_4c_5^3
 {+}c_2c_3^7c_5^2
 {+}c_3^9c_6
 {+}c_3^8c_4c_5
  \\ & \qquad
 {+}c_2^5c_3^2c_5^2c_7
 {+}c_2c_3^8c_7)
  \\ & \quad
 {+}(c_3c_5^5c_6
 {+}c_4c_5^6
 {+}c_2^2c_5^6
 {+}c_2c_3^2c_5^4c_6
 {+}c_2c_3c_4c_5^5
 {+}c_3^3c_5^5
 {+}c_3^2c_4^2c_5^4
 {+}c_2^4c_3c_5c_6^3
  \\ & \qquad
 {+}c_2^4c_4c_5^2c_6^2
 {+}c_2^2c_3^3c_5^3c_6
 {+}c_2^2c_3^2c_4c_5^4
 {+}c_3^5c_4^2c_5c_6
 {+}c_3^4c_4^3c_5^2
 {+}c_2^6c_5^2c_6^2
  \\ & \qquad
 {+}c_2^5c_3^2c_6^3
 {+}c_2^5c_3c_4c_5c_6^2
 {+}c_2^4c_3^3c_5c_6^2
 {+}c_2^4c_3^2c_4^2c_6^2
 {+}c_2^4c_3^2c_5^4
 {+}c_2^3c_3^4c_5^2c_6
  \\ & \qquad
 {+}c_2^3c_3^3c_4c_5^3
 {+}c_2^2c_3^5c_5^3
 {+}c_2c_3^6c_4^2c_6
 {+}c_2c_3^5c_4^3c_5
 {+}c_3^7c_4^2c_5
 {+}c_3^6c_4^4
 {+}c_2c_5^5c_7
  \\ & \qquad
 {+}c_2^2c_3c_5^4c_7
 {+}c_3^4c_5^3c_7
 {+}c_2^5c_5c_6^2c_7
 {+}c_2^4c_3c_5^2c_6c_7
 {+}c_2^4c_4c_5^3c_7
 {+}c_2^3c_3^2c_5^3c_7
  \\ & \qquad
 {+}c_2^2c_3^4c_5c_6c_7
 {+}c_2c_3^4c_4^2c_5c_7
 {+}c_3^6c_4c_5c_7
 {+}c_2^6c_4c_5c_6c_7
 {+}c_2^6c_5^3c_7
  \\ & \qquad
 {+}c_2^5c_3^2c_5c_6c_7
 {+}c_2^5c_3c_4c_5^2c_7
 {+}c_2^4c_3^3c_4c_6c_7
 {+}c_2^4c_3^2c_4^2c_5c_7
 {+}c_2^3c_3^5c_6c_7
  \\ & \qquad
 {+}c_2^2c_3^5c_4^2c_7
 {+}c_2c_3^7c_4c_7
 {+}c_3^9c_7
 {+}c_2^2c_3^2c_5^2c_7^2
 {+}c_3^5c_5c_7^2
 {+}c_2^5c_5^2c_7^2
 {+}c_2^4c_3^2c_6c_7^2
  \\ & \qquad
 {+}c_2^4c_3c_4c_5c_7^2
 {+}c_2c_3^6c_7^2
 {+}c_2^7c_6c_7^2
 {+}c_2^6c_4^2c_7^2
 {+}c_2^5c_3^2c_4c_7^2
 {+}c_2^5c_3c_7^3)
  \\ & \quad
 {+}(c_5^7
 {+}c_2c_3c_5^6
 {+}c_3^3c_5^4c_6
 {+}c_3^2c_4c_5^5
 {+}c_2^4c_5^3c_6^2
 {+}c_2^2c_3^2c_5^5
 {+}c_3^4c_4^2c_5^3
 {+}c_2^5c_3c_5^2c_6^2
  \\ & \qquad
 {+}c_2^4c_3^3c_6^3
 {+}c_2^4c_3^2c_4c_5c_6^2
 {+}c_2^3c_3^3c_5^4
 {+}c_2^2c_3^5c_5^2c_6
 {+}c_2^2c_3^4c_4c_5^3
 {+}c_2c_3^5c_4^2c_5^2
  \\ & \qquad
 {+}c_3^7c_4^2c_6
 {+}c_3^6c_4^3c_5
 {+}c_2c_3^2c_5^4c_7
 {+}c_2^4c_5^4c_7
 {+}c_3^6c_5^2c_7
 {+}c_2^6c_5^2c_6c_7
 {+}c_2^5c_3^2c_6^2c_7
  \\ & \qquad
 {+}c_2^5c_3c_5^3c_7
 {+}c_2^4c_3^2c_4c_5^2c_7
 {+}c_2^3c_3^4c_5^2c_7
 {+}c_2^2c_3^6c_6c_7
 {+}c_2c_3^7c_5c_7
 {+}c_2c_3^6c_4^2c_7
  \\ & \qquad
 {+}c_3^8c_4c_7
 {+}c_2^4c_3c_5^2c_7^2
 {+}c_2^2c_3^4c_5c_7^2
 {+}c_3^7c_7^2
 {+}c_2^6c_3c_6c_7^2
 {+}c_2^6c_4c_5c_7^2
 {+}c_2^3c_3^5c_7^2
  \\ & \qquad
 {+}c_2^7c_7^3)
  , \\
  c(\lambda^6) & =
  1
 {+}c_2
 {+}c_3
 {+}c_4
 {+}c_5
 {+}c_6
 {+}c_7.
\end{align*}
By \eqref{eq:whitney_sum}, we obtain
\begin{align*}
  c_{32}(1 \oplus \lambda^2 \oplus \lambda^4 \oplus \lambda^6) & =
  c_3c_5^3c_7^2
  + c_5^5c_7
  + c_5^4c_6^2
  + c_2^3c_5c_7^3
  + c_2^2c_3c_5c_6c_7^2
  + c_2^2c_4c_5^2c_7^2
  \\ & \quad
  + c_2^2c_5^3c_6c_7
  + c_2c_3^3c_7^3
  + c_3^4c_6c_7^2
  + c_3^3c_4c_5c_7^2
  + c_3^3c_5^2c_6c_7
  \\ & \quad
  + c_2^4c_3c_7^3
  + c_2^4c_5^2c_7^2
  + c_2^4c_6^4
  + c_2^3c_3c_4c_5c_7^2
  + c_2^3c_3c_5^2c_6c_7
  \\ & \quad
  + c_2^2c_3^3c_5c_7^2
  + c_2^2c_3^3c_6^2c_7
  + c_2^2c_3^2c_4^2c_7^2
  + c_2^2c_3^2c_4c_5c_6c_7
  \\ & \quad
  + c_2^2c_3^2c_5^3c_7
  + c_2^2c_3^2c_5^2c_6^2
  + c_2c_3^4c_4c_7^2
  + c_2c_3^4c_5c_6c_7
  + c_3^6c_7^2
  \\ & \quad
  + c_3^5c_4c_6c_7
  + c_3^4c_4^2c_5c_7
  + c_3^4c_4^2c_6^2
  + c_3^4c_5^4
  + c_4^8
  + c_2^{16}
  .
\end{align*}
Recall that $f_{13}^*\Delta_{14}^+ = (\Delta_{13})_{\mathbb C}$, which
gives $w(\smash{f_{13}^*(\Delta_{14}^+)_{\mathbb R})} =
w(\smash{(\Delta_{13})_{\mathbb R}})$.
In order to determine $w_{64}(\smash{(\Delta_{14}^+)_{\mathbb R}})$, we use the
method of indeterminate coefficients;
using the equations $\Sq^nw_{64}(\smash{(\Delta_{14}^+)_{\mathbb R}}) = 0$ for $n = 1, 2, 4$,
we can determine all the coefficients,
and hence the 64--th Stiefel--Whitney class is
given as follows:
\begin{align*}
  w_{64}((\Delta&_{14}^+)_{\mathbb R})  =
  y_{11}y_{13}^3y_{14}
 {+}y_4y_6y_{13}^2y_{14}^2
 {+}y_6^2y_{11}y_{13}y_{14}^2
 {+}y_4y_{10}y_{11}^2y_{14}^2
 {+}y_6y_{10}^3y_{14}^2
  \\ & \ \ 
 {+}y_4y_7y_{13}^3y_{14}
 {+}y_4y_{10}y_{11}y_{12}y_{13}y_{14}
 {+}y_7^2y_{11}y_{12}y_{13}y_{14}
  \\ & \ \ 
 {+}y_7y_8y_{11}^2y_{13}y_{14}
 {+}y_6y_{10}y_{11}^2y_{12}y_{14}
 {+}y_8y_{10}^2y_{11}^2y_{14}
 {+}y_{10}^5y_{14}
  \\ & \ \ 
 {+}y_4^3y_{10}y_{14}^3
 {+}y_4y_6^3y_{14}^3
 {+}y_4^3y_{11}y_{13}y_{14}^2
 {+}y_4^2y_6y_{10}y_{12}y_{14}^2
  \\ & \ \ 
 {+}y_4y_6y_7^2y_{12}y_{14}^2
 {+}y_6^4y_{12}y_{14}^2
 {+}y_4y_6y_7y_8y_{11}y_{14}^2
 {+}y_6^3y_7y_{11}y_{14}^2
  \\ & \ \ 
 {+}y_4^2y_8y_{10}^2y_{14}^2
 {+}y_6^3y_8y_{10}y_{14}^2
 {+}y_4^2y_8^2y_{13}^2y_{14}
 {+}y_4^2y_6y_{11}y_{12}y_{13}y_{14}
  \\ & \ \ 
 {+}y_4y_7^3y_{12}y_{13}y_{14}
 {+}y_4^2y_7y_{11}^2y_{13}y_{14}
 {+}y_4^2y_8y_{10}y_{11}y_{13}y_{14}
  \\ & \ \ 
 {+}y_4y_7^2y_8y_{11}y_{13}y_{14}
 {+}y_7^3y_8^2y_{13}y_{14}
 {+}y_4^2y_7y_{11}y_{12}^2y_{14}
 {+}y_4y_6^2y_{11}^2y_{12}y_{14}
  \\ & \ \ 
 {+}y_6y_7^3y_{11}y_{12}y_{14}
 {+}y_4^2y_{10}^3y_{12}y_{14}
 {+}y_6^3y_{10}^2y_{12}y_{14}
 {+}y_6y_7^2y_8y_{11}^2y_{14}
  \\ & \ \ 
 {+}y_4^4y_6y_{14}^3
 {+}y_4^4y_7y_{13}y_{14}^2
 {+}y_4^3y_6y_7y_{11}y_{14}^2
 {+}y_4^4y_{10}^2y_{14}^2
  \\ & \ \ 
 {+}y_4^3y_6y_8y_{10}y_{14}^2
 {+}y_4^2y_6^3y_{10}y_{14}^2
 {+}y_4^2y_6^2y_8^2y_{14}^2
 {+}y_4y_6^4y_8y_{14}^2
  \\ & \ \ 
 {+}y_4y_6^3y_7^2y_{14}^2
 {+}y_6^6y_{14}^2
 {+}y_4^3y_6y_8y_{11}y_{13}y_{14}
 {+}y_4^3y_7^2y_{11}y_{13}y_{14}
  \\ & \ \ 
 {+}y_4^2y_6y_7y_8^2y_{13}y_{14}
 {+}y_4^2y_6^3y_{12}^2y_{14}
 {+}y_4^3y_6y_{10}^2y_{12}y_{14}
  \\ & \ \ 
 {+}y_4^2y_6^2y_8y_{10}y_{12}y_{14}
 {+}y_4y_6^4y_{10}y_{12}y_{14}
 {+}y_4^2y_7^2y_8^2y_{12}y_{14}
  \\ & \ \ 
 {+}y_4y_6^2y_7^2y_8y_{12}y_{14}
 {+}y_6^5y_8y_{12}y_{14}
 {+}y_4y_6y_7^4y_{12}y_{14}
 {+}y_4^2y_6^2y_8y_{11}^2y_{14}
  \\ & \ \ 
 {+}y_4y_6y_7^3y_8y_{11}y_{14}
 {+}y_4^2y_6^2y_{10}^3y_{14}
 {+}y_6^4y_8^2y_{10}y_{14}
 {+}y_6^2y_7^2y_8^3y_{14}
  \\ & \ \ 
 {+}y_6y_7^4y_8^2y_{14}
 {+}y_{12}y_{13}^4
 {+}y_6^2y_{13}^4
 {+}y_4y_{10}y_{11}y_{13}^3
 {+}y_7^2y_{11}y_{13}^3
  \\ & \ \ 
 {+}y_4y_{11}^2y_{12}y_{13}^2
 {+}y_7y_8y_{11}y_{12}y_{13}^2
 {+}y_6y_{10}^2y_{12}y_{13}^2
 {+}y_6y_{10}y_{11}^2y_{13}^2
  \\ & \ \ 
 {+}y_8y_{10}^2y_{11}y_{12}y_{13}
 {+}y_{10}^4y_{11}y_{13}
 {+}y_4^3y_{13}^4
 {+}y_4^2y_6y_{11}y_{13}^3
 {+}y_4y_7^3y_{13}^3
  \\ & \ \ 
 {+}y_6^3y_8y_{12}y_{13}^2
 {+}y_6^2y_7^2y_{12}y_{13}^2
 {+}y_6^2y_7y_8y_{11}y_{13}^2
 {+}y_6y_7^3y_{11}y_{13}^2
  \\ & \ \ 
 {+}y_4^2y_8y_{11}y_{12}^2y_{13}
 {+}y_4y_6y_7y_{11}^2y_{12}y_{13}
 {+}y_4^2y_{10}^2y_{11}y_{12}y_{13}
  \\ & \ \ 
 {+}y_6^3y_{10}y_{11}y_{12}y_{13}
 {+}y_6^2y_8^2y_{11}y_{12}y_{13}
 {+}y_6y_7^2y_8y_{11}y_{12}y_{13}
  \\ & \ \ 
 {+}y_7^4y_{11}y_{12}y_{13}
 {+}y_7^3y_8y_{11}^2y_{13}
 {+}y_4^3y_6y_7y_{13}^3
 {+}y_4^3y_7^2y_{12}y_{13}^2
  \\ & \ \ 
 {+}y_4^2y_6^3y_{12}y_{13}^2
 {+}y_4^2y_6^2y_7y_{11}y_{13}^2
 {+}y_4^2y_7^2y_8^2y_{13}^2
 {+}y_4y_6^2y_7^2y_8y_{13}^2
  \\ & \ \ 
 {+}y_4y_6y_7^4y_{13}^2
 {+}y_4^3y_6y_{10}y_{11}y_{12}y_{13}
 {+}y_4^2y_6y_7^2y_{11}y_{12}y_{13}
  \\ & \ \ 
 {+}y_4y_6^4y_{11}y_{12}y_{13}
 {+}y_6^4y_7y_8y_{12}y_{13}
 {+}y_4y_7^5y_{12}y_{13}
 {+}y_4^2y_6^2y_{10}^2y_{11}y_{13}
  \\ & \ \ 
 {+}y_6^4y_8^2y_{11}y_{13}
 {+}y_4y_7^4y_8y_{11}y_{13}
 {+}y_7^5y_8^2y_{13}
 {+}y_{10}^4y_{12}^2
 {+}y_{10}^3y_{11}^2y_{12}
  \\ & \ \ 
 {+}y_7^4y_{12}^3
 {+}y_4^2y_{10}y_{11}^2y_{12}^2
 {+}y_4y_7^2y_{11}^2y_{12}^2
 {+}y_6^3y_{11}^2y_{12}^2
 {+}y_7^3y_8y_{11}y_{12}^2
  \\ & \ \ 
 {+}y_7^2y_8^2y_{11}^2y_{12}
 {+}y_4^4y_{12}^4
 {+}y_4^2y_6y_7^2y_{12}^3
 {+}y_4^3y_6y_{11}^2y_{12}^2
 {+}y_4^2y_6y_7y_8y_{11}y_{12}^2
  \\ & \ \ 
 {+}y_4^2y_7^3y_{11}y_{12}^2
 {+}y_4y_6^3y_7y_{11}y_{12}^2
 {+}y_4^2y_6^2y_{10}^2y_{12}^2
 {+}y_6^4y_8^2y_{12}^2
  \\ & \ \ 
 {+}y_6^3y_7^2y_8y_{12}^2
 {+}y_4^2y_6^2y_{10}y_{11}^2y_{12}
 {+}y_4y_7^3y_8^2y_{11}y_{12}
 {+}y_6^3y_7y_8^2y_{11}y_{12}
  \\ & \ \ 
 {+}y_7^4y_8^3y_{12}
 {+}y_6^2y_7^4y_{12}^2
 {+}y_{10}^2y_{11}^4
 {+}y_4^2y_6^2y_{11}^4
 {+}y_6^2y_7^2y_8^2y_{11}^2
 {+}y_7^6y_{11}^2
  \\ & \ \ 
 {+}y_6^4y_{10}^4
 {+}y_8^8
 {+}y_4^{16}
  .
\end{align*}
Thus we obtain the total Stiefel--Whitney class
\begin{align*}
  w((\Delta_{14}^+)_{\mathbb R}) & = c(\Delta_{14}^+) \\
  & =
  1 + w_{64}((\Delta_{14}^+)_{\mathbb R})
  + \Sq^{32}w_{64}((\Delta_{14}^+)_{\mathbb R})
  + \Sq^{48}w_{64}((\Delta_{14}^+)_{\mathbb R})
  \\ & \quad
  + \Sq^{56}w_{64}((\Delta_{14}^+)_{\mathbb R})
  + \Sq^{60}w_{64}((\Delta_{14}^+)_{\mathbb R})
  + \Sq^{62}w_{64}((\Delta_{14}^+)_{\mathbb R})
  \\ & \quad
  + u_{128}.
\end{align*}
Note that $c(\Delta_{14}^+) = c(\Delta_{14}^-)$, since $\Delta_{14}^+$
and $\Delta_{14}^-$ are conjugate to each other.

Recall from \eqref{eq:spin_rep} that
$f_{14}^*\Delta_{15} = (\Delta_{14}^+)_{\mathbb R}$, which gives
$w(f_{14}^*\Delta_{15}) = w((\Delta_{14}^+)_{\mathbb R})$.
In order to determine $w_{64}(\Delta_{15})$, we use the
method of indeterminate coefficients;
using the equations $\Sq^nw_{64}(\Delta_{15}) = 0$ for $n = 1, 2, 4, 8$,
we can determine all the indeterminate coefficients,
and hence the 64--th Stiefel--Whitney class is
given as follows:
\begin{align*}
  w_{64}(\Delta&_{15})  =
  y_{10}y_{13}^3y_{15}
  {+} y_4^2y_{13}^2y_{15}^2
  {+} y_4y_6y_{11}y_{13}y_{15}^2
  {+} y_7^3y_{13}y_{15}^2
  {+} y_4y_{10}^3y_{15}^2
  \\ & \quad
  {+} y_4y_7y_{11}y_{13}y_{14}y_{15}
  {+} y_4y_{10}^2y_{11}y_{14}y_{15}
  {+} y_4y_8y_{11}y_{13}^2y_{15}
  \\ & \quad
  {+} y_6y_7y_{11}y_{12}y_{13}y_{15}
  {+} y_4y_{10}y_{11}^2y_{13}y_{15}
  {+} y_6y_8y_{11}^2y_{13}y_{15}
  {+} y_7^2y_{11}^2y_{13}y_{15}
  \\ & \quad
  {+} y_6y_{10}^3y_{13}y_{15}
  {+} y_6y_{10}^2y_{11}y_{12}y_{15}
  {+} y_8y_{10}^3y_{11}y_{15}
  {+} y_4^3y_7y_{15}^3
  \\ & \quad
  {+} y_4^2y_7^2y_{12}y_{15}^2
  {+} y_4^3y_{11}^2y_{15}^2
  {+} y_4^2y_7y_8y_{11}y_{15}^2
  {+} y_4y_6^2y_7y_{11}y_{15}^2
  {+} y_6y_7^4y_{15}^2
  \\ & \quad
  {+} y_4^2y_6y_7y_{14}^2y_{15}
  {+} y_4^3y_{10}y_{13}y_{14}y_{15}
  {+} y_4y_6^3y_{13}y_{14}y_{15}
  \\ & \quad
  {+} y_4^2y_6y_{10}y_{11}y_{14}y_{15}
  {+} y_4y_6y_7^2y_{11}y_{14}y_{15}
  {+} y_4^2y_6y_{11}^2y_{13}y_{15}
  \\ & \quad
  {+} y_4y_6y_7y_8y_{11}y_{13}y_{15}
  {+} y_4y_6^2y_{10}^2y_{13}y_{15}
  {+} y_4y_6y_8^2y_{10}y_{13}y_{15}
  \\ & \quad
  {+} y_6^3y_8y_{10}y_{13}y_{15}
  {+} y_6y_7^2y_8^2y_{13}y_{15}
  {+} y_7^4y_8y_{13}y_{15}
  {+} y_4^2y_6y_{11}y_{12}^2y_{15}
  \\ & \quad
  {+} y_4y_7^3y_{12}^2y_{15}
  {+} y_6^2y_7^2y_{11}y_{12}y_{15}
  {+} y_4y_7y_8^2y_{11}^2y_{15}
  {+} y_6^2y_7y_8y_{11}^2y_{15}
  \\ & \quad
  {+} y_6y_7^3y_{11}^2y_{15}
  {+} y_4^2y_{10}^3y_{11}y_{15}
  {+} y_6^3y_{10}^2y_{11}y_{15}
  {+} y_4^5y_{14}y_{15}^2
  {+} y_4^4y_8y_{10}y_{15}^2
  \\ & \quad
  {+} y_4^3y_6y_8^2y_{15}^2
  {+} y_4^2y_6^3y_8y_{15}^2
  {+} y_4^2y_6^2y_7^2y_{15}^2
  {+} y_4y_6^5y_{15}^2
  {+} y_4^4y_6y_{13}y_{14}y_{15}
  \\ & \quad
  {+} y_4^4y_7y_{12}y_{14}y_{15}
  {+} y_4^3y_7y_8^2y_{14}y_{15}
  {+} y_4^2y_6^2y_7y_8y_{14}y_{15}
  {+} y_4^4y_7y_{13}^2y_{15}
  \\ & \quad
  {+} y_4^4y_8y_{12}y_{13}y_{15}
  {+} y_4^3y_6^2y_{12}y_{13}y_{15}
  {+} y_4^2y_6^2y_8^2y_{13}y_{15}
  {+} y_6^6y_{13}y_{15}
  \\ & \quad
  {+} y_4^3y_6y_7y_{12}^2y_{15}
  {+} y_4^4y_{10}y_{11}y_{12}y_{15}
  {+} y_4^2y_6y_7y_8^2y_{12}y_{15}
  {+} y_4y_6^2y_7^3y_{12}y_{15}
  \\ & \quad
  {+} y_4^3y_6y_{10}^2y_{11}y_{15}
  {+} y_4y_6^4y_{10}y_{11}y_{15}
  {+} y_4^2y_7^2y_8^2y_{11}y_{15}
  {+} y_4y_6^2y_7^2y_8y_{11}y_{15}
  \\ & \quad
  {+} y_6^5y_8y_{11}y_{15}
  {+} y_4y_6y_7^4y_{11}y_{15}
  {+} y_6^4y_7^2y_{11}y_{15}
  {+} y_4y_7^3y_8^3y_{15}
  {+} y_6^3y_7y_8^3y_{15}
  \\ & \quad
  {+} y_6^2y_7^3y_8^2y_{15}
  {+} y_6y_7^5y_8y_{15}
  {+} y_7^7y_{15}
  {+} y_4^6y_{10}y_{15}^2
  {+} y_4^5y_7^2y_{15}^2
  {+} y_4^4y_6^3y_{15}^2
  \\ & \quad
  {+} y_4^6y_{11}y_{14}y_{15}
  {+} y_4^4y_6^2y_7y_{14}y_{15}
  {+} y_4^6y_{12}y_{13}y_{15}
  {+} y_4^5y_6y_{10}y_{13}y_{15}
  \\ & \quad
  {+} y_4^4y_6^2y_8y_{13}y_{15}
  {+} y_4^4y_6y_7^2y_{13}y_{15}
  {+} y_4^4y_7^3y_{12}y_{15}
  {+} y_4^5y_7y_{11}^2y_{15}
  \\ & \quad
  {+} y_4^4y_7^2y_8y_{11}y_{15}
  {+} y_{11}y_{13}^3y_{14}
  {+} y_4y_6y_{13}^2y_{14}^2
  {+} y_6^2y_{11}y_{13}y_{14}^2
  \\ & \quad
  {+} y_4y_{10}y_{11}^2y_{14}^2
  {+} y_6y_{10}^3y_{14}^2
  {+} y_4y_7y_{13}^3y_{14}
  {+} y_4y_{10}y_{11}y_{12}y_{13}y_{14}
  \\ & \quad
  {+} y_7^2y_{11}y_{12}y_{13}y_{14}
  {+} y_7y_8y_{11}^2y_{13}y_{14}
  {+} y_6y_{10}y_{11}^2y_{12}y_{14}
  {+} y_8y_{10}^2y_{11}^2y_{14}
  \\ & \quad
  {+} y_{10}^5y_{14}
  {+} y_4^3y_{10}y_{14}^3
  {+} y_4y_6^3y_{14}^3
  {+} y_4^3y_{11}y_{13}y_{14}^2
  {+} y_4^2y_6y_{10}y_{12}y_{14}^2
  \\ & \quad
  {+} y_4y_6y_7^2y_{12}y_{14}^2
  {+} y_6^4y_{12}y_{14}^2
  {+} y_4y_6y_7y_8y_{11}y_{14}^2
  {+} y_6^3y_7y_{11}y_{14}^2
  \\ & \quad
  {+} y_4^2y_8y_{10}^2y_{14}^2
  {+} y_6^3y_8y_{10}y_{14}^2
  {+} y_4^2y_8^2y_{13}^2y_{14}
  {+} y_4^2y_6y_{11}y_{12}y_{13}y_{14}
  \\ & \quad
  {+} y_4y_7^3y_{12}y_{13}y_{14}
  {+} y_4^2y_7y_{11}^2y_{13}y_{14}
  {+} y_4^2y_8y_{10}y_{11}y_{13}y_{14}
  \\ & \quad
  {+} y_4y_7^2y_8y_{11}y_{13}y_{14}
  {+} y_7^3y_8^2y_{13}y_{14}
  {+} y_4^2y_7y_{11}y_{12}^2y_{14}
  {+} y_4y_6^2y_{11}^2y_{12}y_{14}
  \\ & \quad
  {+} y_6y_7^3y_{11}y_{12}y_{14}
  {+} y_4^2y_{10}^3y_{12}y_{14}
  {+} y_6^3y_{10}^2y_{12}y_{14}
  {+} y_6y_7^2y_8y_{11}^2y_{14}
  \\ & \quad
  {+} y_4^4y_6y_{14}^3
  {+} y_4^4y_7y_{13}y_{14}^2
  {+} y_4^3y_6y_7y_{11}y_{14}^2
  {+} y_4^4y_{10}^2y_{14}^2
  {+} y_4^3y_6y_8y_{10}y_{14}^2
  \\ & \quad
  {+} y_4^2y_6^3y_{10}y_{14}^2
  {+} y_4^2y_6^2y_8^2y_{14}^2
  {+} y_4y_6^4y_8y_{14}^2
  {+} y_4y_6^3y_7^2y_{14}^2
  {+} y_6^6y_{14}^2
  \\ & \quad
  {+} y_4^3y_6y_8y_{11}y_{13}y_{14}
  {+} y_4^3y_7^2y_{11}y_{13}y_{14}
  {+} y_4^2y_6y_7y_8^2y_{13}y_{14}
  {+} y_4^2y_6^3y_{12}^2y_{14}
  \\ & \quad
  {+} y_4^3y_6y_{10}^2y_{12}y_{14}
  {+} y_4^2y_6^2y_8y_{10}y_{12}y_{14}
  {+} y_4y_6^4y_{10}y_{12}y_{14}
  {+} y_4^2y_7^2y_8^2y_{12}y_{14}
  \\ & \quad
  {+} y_4y_6^2y_7^2y_8y_{12}y_{14}
  {+} y_6^5y_8y_{12}y_{14}
  {+} y_4y_6y_7^4y_{12}y_{14}
  {+} y_4^2y_6^2y_8y_{11}^2y_{14}
  \\ & \quad
  {+} y_4y_6y_7^3y_8y_{11}y_{14}
  {+} y_4^2y_6^2y_{10}^3y_{14}
  {+} y_6^4y_8^2y_{10}y_{14}
  {+} y_6^2y_7^2y_8^3y_{14}
  \\ & \quad
  {+} y_6y_7^4y_8^2y_{14}
  {+} y_{12}y_{13}^4
  {+} y_6^2y_{13}^4
  {+} y_4y_{10}y_{11}y_{13}^3
  {+} y_7^2y_{11}y_{13}^3
  \\ & \quad
  {+} y_4y_{11}^2y_{12}y_{13}^2
  {+} y_7y_8y_{11}y_{12}y_{13}^2
  {+} y_6y_{10}^2y_{12}y_{13}^2
  {+} y_6y_{10}y_{11}^2y_{13}^2
  \\ & \quad
  {+} y_8y_{10}^2y_{11}y_{12}y_{13}
  {+} y_{10}^4y_{11}y_{13}
  {+} y_4^3y_{13}^4
  {+} y_4^2y_6y_{11}y_{13}^3
  {+} y_4y_7^3y_{13}^3
  \\ & \quad
  {+} y_6^3y_8y_{12}y_{13}^2
  {+} y_6^2y_7^2y_{12}y_{13}^2
  {+} y_6^2y_7y_8y_{11}y_{13}^2
  {+} y_6y_7^3y_{11}y_{13}^2
  \\ & \quad
  {+} y_4^2y_8y_{11}y_{12}^2y_{13}
  {+} y_4y_6y_7y_{11}^2y_{12}y_{13}
  {+} y_4^2y_{10}^2y_{11}y_{12}y_{13}
  \\ & \quad
  {+} y_6^3y_{10}y_{11}y_{12}y_{13}
  {+} y_6^2y_8^2y_{11}y_{12}y_{13}
  {+} y_6y_7^2y_8y_{11}y_{12}y_{13}
  {+} y_7^4y_{11}y_{12}y_{13}
  \\ & \quad
  {+} y_7^3y_8y_{11}^2y_{13}
  {+} y_4^3y_6y_7y_{13}^3
  {+} y_4^3y_7^2y_{12}y_{13}^2
  {+} y_4^2y_6^3y_{12}y_{13}^2
  {+} y_4^2y_6^2y_7y_{11}y_{13}^2
  \\ & \quad
  {+} y_4^2y_7^2y_8^2y_{13}^2
  {+} y_4y_6^2y_7^2y_8y_{13}^2
  {+} y_4y_6y_7^4y_{13}^2
  {+} y_4^3y_6y_{10}y_{11}y_{12}y_{13}
  \\ & \quad
  {+} y_4^2y_6y_7^2y_{11}y_{12}y_{13}
  {+} y_4y_6^4y_{11}y_{12}y_{13}
  {+} y_6^4y_7y_8y_{12}y_{13}
  {+} y_4y_7^5y_{12}y_{13}
  \\ & \quad
  {+} y_4^2y_6^2y_{10}^2y_{11}y_{13}
  {+} y_6^4y_8^2y_{11}y_{13}
  {+} y_4y_7^4y_8y_{11}y_{13}
  {+} y_7^5y_8^2y_{13}
  {+} y_{10}^4y_{12}^2
  \\ & \quad
  {+} y_{10}^3y_{11}^2y_{12}
  {+} y_7^4y_{12}^3
  {+} y_4^2y_{10}y_{11}^2y_{12}^2
  {+} y_4y_7^2y_{11}^2y_{12}^2
  {+} y_6^3y_{11}^2y_{12}^2
  \\ & \quad
  {+} y_7^3y_8y_{11}y_{12}^2
  {+} y_7^2y_8^2y_{11}^2y_{12}
  {+} y_4^4y_{12}^4
  {+} y_4^2y_6y_7^2y_{12}^3
  {+} y_4^3y_6y_{11}^2y_{12}^2
  \\ & \quad
  {+} y_4^2y_6y_7y_8y_{11}y_{12}^2
  {+} y_4^2y_7^3y_{11}y_{12}^2
  {+} y_4y_6^3y_7y_{11}y_{12}^2
  {+} y_4^2y_6^2y_{10}^2y_{12}^2
  \\ & \quad
  {+} y_6^4y_8^2y_{12}^2
  {+} y_6^3y_7^2y_8y_{12}^2
  {+} y_4^2y_6^2y_{10}y_{11}^2y_{12}
  {+} y_4y_7^3y_8^2y_{11}y_{12}
  \\ & \quad
  {+} y_6^3y_7y_8^2y_{11}y_{12}
  {+} y_7^4y_8^3y_{12}
  {+} y_6^2y_7^4y_{12}^2
  {+} y_{10}^2y_{11}^4
  {+} y_4^2y_6^2y_{11}^4
  {+} y_6^2y_7^2y_8^2y_{11}^2
  \\ & \quad
  {+} y_7^6y_{11}^2
  {+} y_6^4y_{10}^4
  {+} y_8^8
  {+} y_4^{16}
  .  
\end{align*}
Using the above result on $w_{64}(\Delta_{15})$, we have the
following theorem.
\begin{thm}
  The total Stiefel--Whitney class of the representation
  $\Delta_{15} \co \Spin(15) \to O(128)$ is given by
  \begin{align*}
    w(\Delta_{15}) & =
    1
    + w_{64}(\Delta_{15})
    + \Sq^{32}w_{64}(\Delta_{15})
    + \Sq^{48}w_{64}(\Delta_{15})
    + \Sq^{56}w_{64}(\Delta_{15})
    \\ & \quad
    + \Sq^{60}w_{64}(\Delta_{15})
    + \Sq^{62}w_{64}(\Delta_{15})
    + \Sq^{63}w_{64}(\Delta_{15})
    + u_{128},
  \end{align*}
  where $u_{128} = w_{128}(\Delta_{15})$.
\end{thm}

\section[The Stiefel--Whitney classes of the induced representation
  from the adjoint representation of E8]{The Stiefel--Whitney classes of the induced representation
  from the adjoint representation of $E_8$}
\label{sec:adjoint_rep}

Summing up all the calculations in the previous section,
we have the following:
\begin{thm}
  \label{thm:adjoint}
    For the
      representation $\Delta_{15} \oplus \lambda_{15}^1 \oplus \lambda_{15}^2
        \co \Spin(15) \to SO(248)$, the Stiefel--Whitney classes of degree $2^i$ are given as follows:
\begin{align*}
  w_1(\Delta_{15} \oplus &\lambda_{15}^1 \oplus \lambda_{15}^2)  =
  w_1(\lambda_{15}^1 \oplus \lambda_{15}^2) =
  0, \hspace{2in}\\
  w_2(\Delta_{15} \oplus &\lambda_{15}^1 \oplus \lambda_{15}^2) =
  w_2(\lambda_{15}^1 \oplus \lambda_{15}^2) =
  0, \\
  w_4(\Delta_{15} \oplus &\lambda_{15}^1 \oplus \lambda_{15}^2)  =
  w_4(\lambda_{15}^1 \oplus \lambda_{15}^2) =
  0, \\
  w_8(\Delta_{15} \oplus &\lambda_{15}^1 \oplus \lambda_{15}^2)  =
  w_8(\lambda_{15}^1 \oplus \lambda_{15}^2) =
  0, \\
  w_{16}(\Delta_{15} \oplus &\lambda_{15}^1 \oplus \lambda_{15}^2)  =
  w_{16}(\lambda_{15}^1 \oplus \lambda_{15}^2) =
  y_4^4
  , \\
  w_{32}(\Delta_{15} \oplus &\lambda_{15}^1 \oplus \lambda_{15}^2)  =
  w_{32}(\lambda_{15}^1 \oplus \lambda_{15}^2)
  \\ & =
  y_4^2y_{11}y_{13}
  + y_6^2y_7y_{13}
  + y_6^2y_{10}^2
  + y_7^3y_{11}
  + y_4^8
  , \\
  w_{64}(\Delta_{15} \oplus &\lambda_{15}^1 \oplus \lambda_{15}^2)  =
  w_{64}(\Delta_{15}) + w_{64}(\lambda_{15}^1 \oplus \lambda_{15}^2)
\\
 & =
  y_{10}y_{13}^3y_{15}
  + y_{11}y_{13}^3y_{14}
  + y_{12}y_{13}^4
  + y_4^2y_{13}^2y_{15}^2
  \\ & \quad
  + y_4y_6y_{11}y_{13}y_{15}^2
  + y_4y_6y_{13}^2y_{14}^2
  + y_4y_7y_{11}y_{13}y_{14}y_{15}
  \\ & \quad
  + y_4y_7y_{13}^3y_{14}
  + y_4y_8y_{11}y_{13}^2y_{15}
  + y_4y_{10}^3y_{15}^2
  \\ & \quad
  + y_4y_{10}^2y_{11}y_{14}y_{15}
  + y_4y_{10}y_{11}^2y_{13}y_{15}
  + y_4y_{10}y_{11}^2y_{14}^2
  \\ & \quad
  + y_4y_{10}y_{11}y_{12}y_{13}y_{14}
  + y_4y_{11}^2y_{12}y_{13}^2
  + y_6^2y_{11}y_{13}y_{14}^2
  \\ & \quad
  + y_6y_7y_{11}y_{12}y_{13}y_{15}
  + y_6y_8y_{11}^2y_{13}y_{15}
  + y_6y_{10}^3y_{13}y_{15}
  \\ & \quad
  + y_6y_{10}^3y_{14}^2
  + y_6y_{10}^2y_{11}y_{12}y_{15}
  + y_6y_{10}^2y_{12}y_{13}^2
  \\ & \quad
  + y_6y_{10}y_{11}^2y_{12}y_{14}
  + y_7^3y_{13}y_{15}^2
  + y_7^2y_{11}^2y_{13}y_{15}
  \\ & \quad
  + y_7^2y_{11}y_{12}y_{13}y_{14}
  + y_7y_8y_{11}^2y_{13}y_{14}
  + y_7y_8y_{11}y_{12}y_{13}^2
  \\ & \quad
  + y_8y_{10}^3y_{11}y_{15}
  + y_8y_{10}^2y_{11}^2y_{14}
  + y_8y_{10}^2y_{11}y_{12}y_{13}
  + y_{10}^5y_{14}
  \\ & \quad
  + y_{10}^4y_{12}^2
  + y_{10}^3y_{11}^2y_{12}
  + y_4^3y_7y_{15}^3
  + y_4^3y_{10}y_{13}y_{14}y_{15}
  \\ & \quad
  + y_4^3y_{10}y_{14}^3
  + y_4^3y_{11}^2y_{15}^2
  + y_4^3y_{11}y_{13}^2y_{15}
  + y_4^3y_{11}y_{13}y_{14}^2
  \\ & \quad
  + y_4^2y_6y_7y_{14}^2y_{15}
  + y_4^2y_6y_{10}y_{11}y_{14}y_{15}
  + y_4^2y_6y_{10}y_{12}y_{14}^2
  \\ & \quad
  + y_4^2y_6y_{11}y_{12}^2y_{15}
  + y_4^2y_6y_{11}y_{12}y_{13}y_{14}
  + y_4^2y_7^2y_{12}y_{15}^2
  \\ & \quad
  + y_4^2y_7y_8y_{11}y_{15}^2
  + y_4^2y_7y_{11}y_{12}^2y_{14}
  + y_4^2y_7y_{11}y_{12}y_{13}^2
  \\ & \quad
  + y_4^2y_8^2y_{13}^2y_{14}
  + y_4^2y_8y_{10}^2y_{14}^2
  + y_4^2y_8y_{10}y_{11}y_{13}y_{14}
  \\ & \quad
  + y_4^2y_8y_{11}y_{12}^2y_{13}
  + y_4^2y_{10}^3y_{11}y_{15}
  + y_4^2y_{10}^3y_{12}y_{14}
  \\ & \quad
  + y_4^2y_{10}^2y_{11}y_{12}y_{13}
  + y_4^2y_{10}y_{11}^2y_{12}^2
  + y_4y_6^3y_{13}y_{14}y_{15}
  \\ & \quad
  + y_4y_6^3y_{14}^3
  + y_4y_6^2y_7y_{11}y_{15}^2
  + y_4y_6^2y_{10}^2y_{13}y_{15}
  \\ & \quad
  + y_4y_6^2y_{11}^2y_{12}y_{14}
  + y_4y_6y_7^2y_{11}y_{14}y_{15}
  + y_4y_6y_7^2y_{12}y_{14}^2
  \\ & \quad
  + y_4y_6y_7y_8y_{11}y_{13}y_{15}
  + y_4y_6y_7y_8y_{11}y_{14}^2
  \\ & \quad
  + y_4y_6y_7y_{11}^2y_{12}y_{13}
  + y_4y_6y_8^2y_{10}y_{13}y_{15}
  + y_4y_7^3y_{11}y_{13}y_{15}
  \\ & \quad
  + y_4y_7^3y_{12}^2y_{15}
  + y_4y_7^3y_{12}y_{13}y_{14}
  + y_4y_7^2y_8y_{11}y_{13}y_{14}
  \\ & \quad
  + y_4y_7^2y_{11}^2y_{12}^2
  + y_4y_7y_8^2y_{11}^2y_{15}
  + y_6^4y_{12}y_{14}^2
  \\ & \quad
  + y_6^3y_7y_{11}y_{13}y_{15}
  + y_6^3y_7y_{11}y_{14}^2
  + y_6^3y_8y_{10}y_{13}y_{15}
  \\ & \quad
  + y_6^3y_8y_{10}y_{14}^2
  + y_6^3y_8y_{12}y_{13}^2
  + y_6^3y_{10}^2y_{11}y_{15}
  \\ & \quad
  + y_6^3y_{10}^2y_{12}y_{14}
  + y_6^3y_{10}y_{11}y_{12}y_{13}
  + y_6^3y_{11}^2y_{12}^2
  \\ & \quad
  + y_6^2y_7^2y_{11}y_{12}y_{15}
  + y_6^2y_7^2y_{11}y_{13}y_{14}
  + y_6^2y_7^2y_{12}y_{13}^2
  \\ & \quad
  + y_6^2y_7y_8y_{11}^2y_{15}
  + y_6^2y_8^2y_{11}y_{12}y_{13}
  + y_6y_7^4y_{15}^2
  \\ & \quad
  + y_6y_7^3y_{11}^2y_{15}
  + y_6y_7^3y_{11}y_{12}y_{14}
  + y_6y_7^2y_8^2y_{13}y_{15}
  \\ & \quad
  + y_6y_7^2y_8y_{11}^2y_{14}
  + y_6y_7^2y_8y_{11}y_{12}y_{13}
  + y_7^4y_8y_{13}y_{15}
  \\ & \quad
  + y_7^4y_{12}^3
  + y_7^3y_8^2y_{13}y_{14}
  + y_7^3y_8y_{11}y_{12}^2
  + y_7^2y_8^2y_{11}^2y_{12}
  \\ & \quad
  + y_4^5y_{14}y_{15}^2
  + y_4^4y_6y_{13}y_{14}y_{15}
  + y_4^4y_6y_{14}^3
  + y_4^4y_7y_{11}y_{15}^2
  \\ & \quad
  + y_4^4y_7y_{12}y_{14}y_{15}
  + y_4^4y_8y_{10}y_{15}^2
  + y_4^4y_8y_{12}y_{13}y_{15}
  \\ & \quad
  + y_4^4y_{10}y_{11}y_{12}y_{15}
  + y_4^4y_{11}y_{12}^2y_{13}
  + y_4^3y_6^2y_{12}y_{13}y_{15}
  \\ & \quad
  + y_4^3y_6y_7y_{11}y_{13}y_{15}
  + y_4^3y_6y_7y_{11}y_{14}^2
  + y_4^3y_6y_7y_{12}^2y_{15}
  \\ & \quad
  + y_4^3y_6y_8^2y_{15}^2
  + y_4^3y_6y_8y_{10}y_{14}^2
  + y_4^3y_6y_8y_{11}y_{13}y_{14}
  \\ & \quad
  + y_4^3y_6y_{10}^2y_{11}y_{15}
  + y_4^3y_6y_{10}^2y_{12}y_{14}
  + y_4^3y_6y_{10}y_{11}y_{12}y_{13}
  \\ & \quad
  + y_4^3y_6y_{11}^2y_{12}^2
  + y_4^3y_7y_8^2y_{14}y_{15}
  + y_4^2y_6^3y_8y_{15}^2
  + y_4^2y_6^3y_{10}y_{14}^2
  \\ & \quad
  + y_4^2y_6^3y_{12}^2y_{14}
  + y_4^2y_6^3y_{12}y_{13}^2
  + y_4^2y_6^2y_7y_8y_{14}y_{15}
  \\ & \quad
  + y_4^2y_6^2y_8^2y_{13}y_{15}
  + y_4^2y_6^2y_8^2y_{14}^2
  + y_4^2y_6^2y_8y_{10}y_{12}y_{14}
  \\ & \quad
  + y_4^2y_6^2y_8y_{11}^2y_{14}
  + y_4^2y_6^2y_{10}^3y_{14}
  + y_4^2y_6^2y_{10}^2y_{12}^2
  \\ & \quad
  + y_4^2y_6^2y_{10}y_{11}^2y_{12}
  + y_4^2y_6y_7^2y_{11}y_{12}y_{13}
  + y_4^2y_6y_7^2y_{12}^3
  \\ & \quad
  + y_4^2y_6y_7y_8^2y_{12}y_{15}
  + y_4^2y_6y_7y_8^2y_{13}y_{14}
  + y_4^2y_6y_7y_8y_{11}y_{12}^2
  \\ & \quad
  + y_4^2y_7^3y_{11}y_{12}^2
  + y_4^2y_7^2y_8^2y_{11}y_{15}
  + y_4^2y_7^2y_8^2y_{12}y_{14}
  \\ & \quad
  + y_4^2y_8^4y_{11}y_{13}
  + y_4y_6^5y_{15}^2
  + y_4y_6^4y_8y_{14}^2
  + y_4y_6^4y_{10}y_{11}y_{15}
  \\ & \quad
  + y_4y_6^4y_{10}y_{12}y_{14}
  + y_4y_6^4y_{11}y_{12}y_{13}
  + y_4y_6^3y_7^2y_{13}y_{15}
  \\ & \quad
  + y_4y_6^3y_7^2y_{14}^2
  + y_4y_6^3y_7y_{11}y_{12}^2
  + y_4y_6^2y_7^3y_{12}y_{15}
  \\ & \quad
  + y_4y_6^2y_7^3y_{13}y_{14}
  + y_4y_6^2y_7^2y_8y_{11}y_{15}
  + y_4y_6^2y_7^2y_8y_{12}y_{14}
  \\ & \quad
  + y_4y_6y_7^4y_{11}y_{15}
  + y_4y_6y_7^4y_{12}y_{14}
  + y_4y_6y_7^3y_8y_{11}y_{14}
  \\ & \quad
  + y_4y_7^3y_8^3y_{15}
  + y_4y_7^3y_8^2y_{11}y_{12}
  + y_6^5y_8y_{11}y_{15}
  + y_6^5y_8y_{12}y_{14}
  \\ & \quad
  + y_6^4y_7^2y_{13}^2
  + y_6^4y_7y_8y_{12}y_{13}
  + y_6^4y_8^2y_{10}y_{14}
  + y_6^4y_8^2y_{12}^2
  \\ & \quad
  + y_6^4y_{10}^4
  + y_6^3y_7^2y_8y_{12}^2
  + y_6^3y_7y_8^3y_{15}
  + y_6^3y_7y_8^2y_{11}y_{12}
  \\ & \quad
  + y_6^2y_7^3y_8^2y_{15}
  + y_6^2y_7^2y_8^3y_{14}
  + y_6^2y_7y_8^4y_{13}
  + y_6^2y_8^4y_{10}^2
  \\ & \quad
  + y_6y_7^5y_8y_{15}
  + y_6y_7^4y_8^2y_{14}
  + y_7^6y_{11}^2
  + y_7^4y_8^3y_{12}
  + y_7^3y_8^4y_{11}
  \\ & \quad
  + y_4^5y_7y_{11}^2y_{15}
  + y_4^5y_{10}^2y_{11}y_{13}
  + y_4^5y_{11}^4
  + y_4^4y_6y_7y_{10}^2y_{15}
  \\ & \quad
  + y_4^4y_6y_7y_{11}^2y_{13}
  + y_4^4y_7^2y_{10}^2y_{14}
  + y_4^4y_7^2y_{11}^2y_{12}
  \\ & \quad
  + y_4^4y_7y_8y_{10}^2y_{13}
  + y_4^4y_7y_8y_{11}^3
  + y_4^7y_6y_{15}^2
  + y_4^7y_7y_{14}y_{15}
  \\ & \quad
  + y_4^7y_{10}y_{11}y_{15}
  + y_4^7y_{10}y_{13}^2
  + y_4^7y_{11}^2y_{14}
  + y_4^7y_{11}y_{12}y_{13}
  \\ & \quad
  + y_4^6y_6^2y_{13}y_{15}
  + y_4^6y_6^2y_{14}^2
  + y_4^6y_6y_7y_{12}y_{15}
  + y_4^6y_6y_7y_{13}y_{14}
  \\ & \quad
  + y_4^6y_6y_8y_{13}^2
  + y_4^6y_7^2y_{12}y_{14}
  + y_4^6y_7^2y_{13}^2
  + y_4^6y_7y_8y_{10}y_{15}
  \\ & \quad
  + y_4^6y_7y_8y_{11}y_{14}
  + y_4^6y_7y_8y_{12}y_{13}
  + y_4^6y_8^2y_{11}y_{13}
  + y_4^6y_{10}^4
  \\ & \quad
  + y_4^4y_6^4y_{11}y_{13}
  + y_4^4y_6^4y_{12}^2
  + y_4^4y_6^2y_8^2y_{10}^2
  + y_4^4y_7^5y_{13}
  \\ & \quad
  + y_4^{10}y_{11}y_{13}
  + y_4^8y_6^2y_7y_{13}
  + y_4^8y_6^2y_{10}^2
  + y_4^8y_7^3y_{11}
  + y_4^4y_6^8
  \\ & \quad
  + y_4^{16}
  , \\
  w_{128}(\Delta_{15} &\oplus \lambda_{15}^1 \oplus \lambda_{15}^2)  =
  w_{128}(\Delta_{15})
  + w_{112}(\Delta_{15})w_{16}(\lambda_{15}^1 \oplus \lambda_{15}^2)
  \\ & \quad
  + w_{96}(\Delta_{15})w_{32}(\lambda_{15}^1 \oplus \lambda_{15}^2)
  + w_{64}(\Delta_{15})w_{64}(\lambda_{15}^1 \oplus \lambda_{15}^2)
  \\ & \hphantom{\, \oplus \lambda_{15}^1 \oplus \lambda_{15}^2)\ }  \equiv
  u_{128}
  \mod {\rm decomposables.}
\end{align*}
\end{thm}

Since the element $w_{128}(\Delta_{15} \oplus \lambda_{15}^1 \oplus
\lambda_{15}^2)$
is a member of a system of generators of
$H^*(B\Spin(15);\mathbb Z/2)$ as an algebra over the Steenrod algebra,
we obtain the following corollary.

\begin{cor}
  \label{cor:gen}
  The Stiefel--Whitney class $w_{128}(\mathrm{Ad}_{E_8})$ of the
  adjoint representation $\mathrm{Ad}_{E_8} \co E_8 \to SO(248)$ can
  be chosen as a member of a system of generators of
  $H^*(BE_8;\mathbb Z/2)$ as an algebra over the Steenrod algebra.
\end{cor}

\bibliographystyle{gtart}
\bibliography{link}

\begin{thebibliography}{}
\providecommand\bibmarginpar{\leavevmode\marginpar}
\def\urlstyle#1{{\tt #1}}

\bibitem{adams1}
\textbf{J\,F Adams}, \emph{Lectures on {L}ie groups}, W. A. Benjamin,, New
  York-Amsterdam (1969) \xox{MR}{0252560}

\bibitem{adams2}
\textbf{J\,F Adams}, \emph{Lectures on exceptional {L}ie groups}, Chicago
  Lectures in Mathematics, University of Chicago Press, Chicago, IL (1996)
  \xox{MR}{1428422}With a foreword by J. Peter May, Edited by Zafer Mahmud and
  Mamoru Mimura

\bibitem{ABS}
\textbf{M\,F Atiyah}, \textbf{R Bott}, \textbf{A Shapiro},
  \href{http://dx.doi.org/10.1016/0040-9383(64)90003-5} {\emph{Clifford
  modules}}, Topology 3 (1964) 3--38 \xox{MR}{0167985}

\bibitem{borel}
\textbf{A Borel},
  \href{http://links.jstor.org/sici?sici=0002-9327(195404)76:2%3C273:SLELCD%3E%
2.0.CO%3B2-J} {\emph{Sur l'homologie et la cohomologie des groupes de {L}ie
  compacts connexes}}, Amer. J. Math. 76 (1954) 273--342 \xox{MR}{0064056}

\bibitem{BH1}
\textbf{A Borel}, \textbf{F Hirzebruch},
  \href{http://links.jstor.org/sici?sici=0002-9327(195804)80:2%3C458:CCAHSI%3E%
2.0.CO%3B2-0} {\emph{Characteristic classes and homogeneous spaces. {I}}},
  Amer. J. Math. 80 (1958) 458--538 \xox{MR}{0102800}

\bibitem{CLO}
\textbf{D Cox}, \textbf{J Little}, , \textbf{D O'Shea}, \emph{Ideals,
  varieties, and algorithms}, second edition, Undergraduate Texts in
  Mathematics, Springer, New York (1997) \xox{MR}{1417938}An introduction to
  computational algebraic geometry and commutative algebra

\bibitem{KM}
\textbf{A Kono}, \textbf{M Mimura},
  \href{http://dx.doi.org/10.1016/0022-4049(75)90013-4} {\emph{Cohomology
  {${\rm mod}\,2$} of the classifying space of th compact connected {L}ie group
  of type {$E\sb{6}$}}}, J. Pure Appl. Algebra 6 (1975) 61--81
  \xox{MR}{0368002}

\bibitem{KMS2}
\textbf{A Kono}, \textbf{M Mimura}, \textbf{N Shimada},
  \href{http://dx.doi.org/10.1016/0022-4049(76)90048-7} {\emph{On the
  cohomology mod {$2$} of the classifying space of the {$1$}\ connected
  exceptional {L}ie group {$E\sb{7}$}}}, J. Pure Appl. Algebra 8 (1976)
  267--283 \xox{MR}{0415652}

\bibitem{Mori}
\textbf{M Mori}, \emph{Computation of $\mathrm{Cotor}_A(\mathbb Z/2, \mathbb
  Z/2)$ for $A = H^*(E_8;\mathbb Z/2)$}, hand written manuscript in Japanese

\bibitem{quillen}
\textbf{D Quillen}, \href{http://dx.doi.org/10.1007/BF01350050} {\emph{The
  {${\rm mod}$} {$2$} cohomology rings of extra-special {$2$}-groups and the
  spinor groups}}, Math. Ann. 194 (1971) 197--212 \xox{MR}{0290401}

\bibitem{toda}
\textbf{H Toda}, \emph{Cohomology of the classifying space of exceptional {L}ie
  groups}, from: ``Manifolds---Tokyo 1973 (Proc. Internat. Conf., Tokyo,
  1973)'', Univ. Tokyo Press, Tokyo (1975)  265--271 \xox{MR}{0368059}

\end{thebibliography}

\end{document}